\documentclass[final,3p,times]{elsarticle}
\usepackage[toc,page,title,titletoc,header]{appendix}
\usepackage{amsmath}
\usepackage{amssymb}
\usepackage{amsthm}
\usepackage{enumerate}
\usepackage{enumitem}
\usepackage{mathrsfs}
\usepackage{graphicx}
\usepackage{epsfig}
\usepackage{tikz,pgfplots,tikz-3dplot}
\usepackage{epstopdf}
\usetikzlibrary{fit}
\biboptions{sort&compress}

\newtheorem{thm}{Theorem}[section]
\newtheorem{lemma}[thm]{Lemma}
\newtheorem{rem}{Remark}

\newcommand{\bR}{{\mathbb R}}
\newcommand{\bZ}{\mathbf{Z}}

\usepackage{mathptmx}
\usepackage{newtxtext}
\usepackage{newtxmath}

\DeclareMathAlphabet{\mathcal}{OMS}{cmsy}{m}{n}

\newcommand{\bI}{\mathbb{I}}
\newcommand{\rd}{\;{\rm d}}
\newcommand{\ds}{\rd s}

\newcommand{\dd}[1]{\frac{\rm d}{{\rm d}#1}}
\newcommand{\ddt}{\dd{t}}
\newcommand{\nn}{\nonumber}
\newcommand{\ttau}{\Delta t}
\newcommand{\mV}{\mathcal{V}}
\newcommand{\hah}{\hat{\vec h}}
\newcommand{\han}{\hat{\vec n}}

\def\epsilon{\varepsilon} 
\renewcommand{\vec}[1]{{\boldsymbol{#1}}}

\begin{document}
\begin{frontmatter}
\title{%\texttt{\jobname} \textbf{DRAFT} \quad {\rm \mydate} \\
A structure-preserving parametric finite element method for solid-state 
dewetting on curved substrates
}

\author[1]{Weizhu Bao}
\address[1]{Department of Mathematics, National University of Singapore, Singapore, 119076}
\ead{matbaowz@nus.edu.sg}
\author[1]{Yifei Li}
\ead{e0444158@u.nus.edu}
\author[2]{Quan Zhao}
\address[2]{School of Mathematical Sciences, University of Science and Technology of China, Hefei 230026, China}
\ead{quanzhao@ustc.edu.cn}

\begin{abstract}
We consider a two-dimensional sharp-interface model for solid-state dewetting of thin films with anisotropic surface energies on curved substrates, where the film/vapor interface and substrate surface are represented by an evolving and a static curve, respectively. The model is governed by the anisotropic surface diffusion for the evolving curve,  with appropriate boundary conditions at the contact points where the two curves meet.  The continuum model obeys an energy decay law and preserves the enclosed area between the two curves. We introduce an arclength parameterization for the substrate curve, which plays a crucial role in a structure-preserving approximation as it straightens the curved substrate and tracks length changes between contact points. Based on this insight, we introduce a symmetrized weak formulation which leads to an unconditional energy stable parametric approximation in terms of the discrete energy. We also provide an error estimate of the enclosed area, which depends on the substrate profile and can be zero in the case of a flat substrate. Furthermore, we introduce a correction to the discrete normals to enable an exact area preservation for general curved substrates. The resulting nonlinear system is efficiently solved using a hybrid iterative algorithm which combines both Picard and Newton's methods. Numerical results are presented to show the robustness and good properties of the introduced method for simulating  solid-state dewetting on various curved substrates. 
\end{abstract} 


\end{frontmatter}

\setcounter{equation}{0}
\section{Introduction} \label{sec:intro}

Deposited solid thin films are typically unstable and could dewet or agglomerate to form isolated islands when they are heated well below the melting points. This phenomenon is called solid-state dewetting (SSD) since the thin films remain in the solid state during the process \cite{Thompson12solid, Leroy16}. In recent decades, SSD has found a various technological applications in optical and magnetic devices, thin films,  sensors and catalyst formation, etc \cite{Armelao06, Schmidt09silicon, Bollani19templated}. A lot of efforts has been devoted to SSD both in the experimental aspect \cite{Jiran90, Jiran92,Giermann05,Cheng06templated, Wang11, Wang2013solid, Lu16nanostructure, Kovalenko17,Naffouti16} and in mathematical modeling and simulation \cite{Srolovitz86b, Dornel06, jiang2012phase, Zucker13, Jiang18curved, Jiang19a, JiangZB20, Boc22stress, Garcke23diffuse}, etc. In particular, numerical simulations play a crucial role in understanding the dewetting patterns of thin films, and thus predictions \cite{jiang2012phase, Backofen2017convexity, Zhao2020p,Boc22stress}.

Experimental results reveals that the morphology changes in SSD are mainly due to the mass transport via surface diffusion flow. In the isotropic materials, the normal velocity for the interface that separates the thin film and the surrounding vapor is proportional to the surface Laplacian of the local mean curvature \cite{Mullins57}. SSD also involves the motion of the contact line, which is formed when the evolving interface meets the substrate. In particular, the equilibrium contact angle at the contact line is given by the Young's law. The first modeling effort for SSD can be tracked back to the work by Sorlovitz and Safran \cite{Srolovitz86b}, where an axisymmetric sharp-interface model was introduced to study the hole growth. The diffuse-interface approach was employed for modelling SSD in \cite{jiang2012phase, Backofen2017convexity}. Later efforts aim to include the anisotropic effects via the discrete model \cite{Dornel06}, the crystalline model \cite{Zucker13},  the diffuse-interface models \cite{Dziwnik17anisotropic, Garcke23diffuse}, and the sharp-interface models \cite{Jiang19a, Jiang18curved, JiangZB20, Boc22stress}. Nevertheless, current theoretical study mostly focus on the case of a flat substrate, and there haven't be well studied for the case of a topologically patterned substrate compared to the experimental results \cite{ahn80, Cheng06templated, Wang11, Wang2013solid, Lu16nanostructure, Kovalenko17}. In the sharp-interface model of SSD,  the evolving interface is represented by an open curve/surface with moving boundaries attached to a fixed substrate surface. The motion for the film/vapor interface is described by (anisotropic) surface diffusion flow with appropriate boundary conditions at the moving contact lines.

The framework of parametric approach has been considered as one of the prominent methods for solving the geometric evolution equations, and the readers can refer to \cite{Bansch05, Deckelnick05, Barrett07, Hausser07, Barrett08JCP, Pozzi08, Kovacs2020, Barrett20, Jiang2021, Zhao2021, BJY21symmetrized, Duan2023} and references therein. As the $H^{-1}$ gradient flow of the surface energy,  the surface diffusion flow  obeys two geometric laws, i.e., the energy dissipation and the mass conservation. Therefore it is of great interest to design numerical approximations so that the two geometric laws are satisfied as well.  This was firstly made possible with an unconditionally stable parametric method in \cite{Bansch05}. Later Barrett, Garcke and N\"urnberg introduced a novel  weak formulation which allows tangential degrees of freedom. This leads to a parametric method which not only satisfies an unconditional stability estimate but also has good mesh properties. The `BGN' method was then generalized to the anisotropic case in the context of anisotropy in the Riemannian metric form \cite{Barrett07Ani,Barrett08Ani}, and a discrete energy decay law was achieved as well.  Based on the  `BGN' formulation and a suitable stabilizing function, Bao et al., \cite{BJY21symmetrized, BL22} introduced a symmetrized weak formulation for anisotropic surface diffusion. This then leads to a parametric finite element method, where an unconditional stability can be established  under more general anisotropies. We also refer the reader to a recent structure-preserving method by Bao and Zhao in \cite{Zhao2021}. In this work, the authors employed a time-weighted discrete normals to achieve an exact mass preservation on fully discrete level.

Parametric approximations have also been employed to open curve/surface evolution, see   \cite{Deckelnick98finite, Barrett07b, Barrett10, garcke2010parametric,BGNZ23} with static contact angle conditions and \cite{ Bao17, Zhao2020p,Zhao20, Bao23energy} with relaxed contact angle conditions. In particular, for the flat boundaries, the introduced methods \cite{Zhao20, Bao23energy, BGNZ23} can also be shown to satisfy an unconditional stability of the total energy including both the film/vapor interface energy and the substrate energy. However, in the case of curved boundaries, the problem becomes much more complicated. When the contact angle condition is weakly enforced in the `BGN' formulation, the velocity of the contact points/lines needs to belong to the tangential space of the substrate surface in order to satisfy the attachment condition, see \eqref{eq:bd1}. On fully discrete level, the discrete vertices on the contact points/lines may not lie on the curved substrate exactly at later time steps. A possible remedy is to employ an orthogonal projection onto the substrate surface at every time step, which however will destroy the stability estimate and mass preservation, see \cite[Remark 3.5]{BGNZ23}. Following the work in \cite{Bao17}, an explicit update of the contact points was employed in \cite{Jiang18curved} for SSD on a curved substrate. This method ignored the variational structure of the considered problem and thus leads to an unnatural treatment with the motion of the contact points. The extension of the parametric method to curved substrates remains an open challenge.

To address the evolution of contact points on curved substrates, we propose a novel approach based on two key insights. Firstly, we introduce an appropriate map to \textit{straighten} the curved substrate into a flat one. This allows for direct updates of the contact points while still ensuring the attachment condition, thus avoiding the undesired orthogonal projections that would destroy the intrinsic variational structure. Secondly, we  adopt the \textit{arclength parameterization} as the specific map for the curved substrate. This parameterization allows to track the change of the length between the contact points, and thus preserves the dissipation of energy along the substrate. These two insights are crucial as they essentially transform the problem of curved substrates into one of flat substrates.

Building upon this understanding, in the current work, we will consider a two-dimensional (2D) sharp-interface model for SSD with anisotropic surface energies on general curved substrates \cite{Jiang18curved}. We combine the idea in \cite{BJY21symmetrized} and \cite{Barrett07b} to propose a weak formulation for the model, which could lead to an unconditionally stable approximation in term of the total free energy. Additionally, we provide an error estimate for mass preservation, which in the case of a flat substrate can be shown to be zero \cite{Zhao2021}. Besides, we show that a suitable adaption of the discrete normals will enable us to  achieve an exact mass preservation.   

The rest of the paper is organized as follows. In Section \ref{sec:model}, we present the sharp-interface model for SSD on general curved substrates. We then propose a weak formulation for the model and prove the energy decay and mass preservation laws within the weak formulation in Section \ref{sec:weak}. Next in Section \ref{sec:fem}, we present a structure-preserving parametric finite element method for the weak formulation, where the energy decay and mass preservation are satisfied well on the discrete level. Subsequently, numerical results are presented to show the robustness of the introduced method and their good properties in Section \ref{sec:num}. Finally, we draw some conclusions in Section \ref{sec:con}. 

\setcounter{equation}{0}
\section{The sharp-interface model}
\label{sec:model}

\begin{figure}[!htp]
\center
\includegraphics[width=0.8\textwidth]{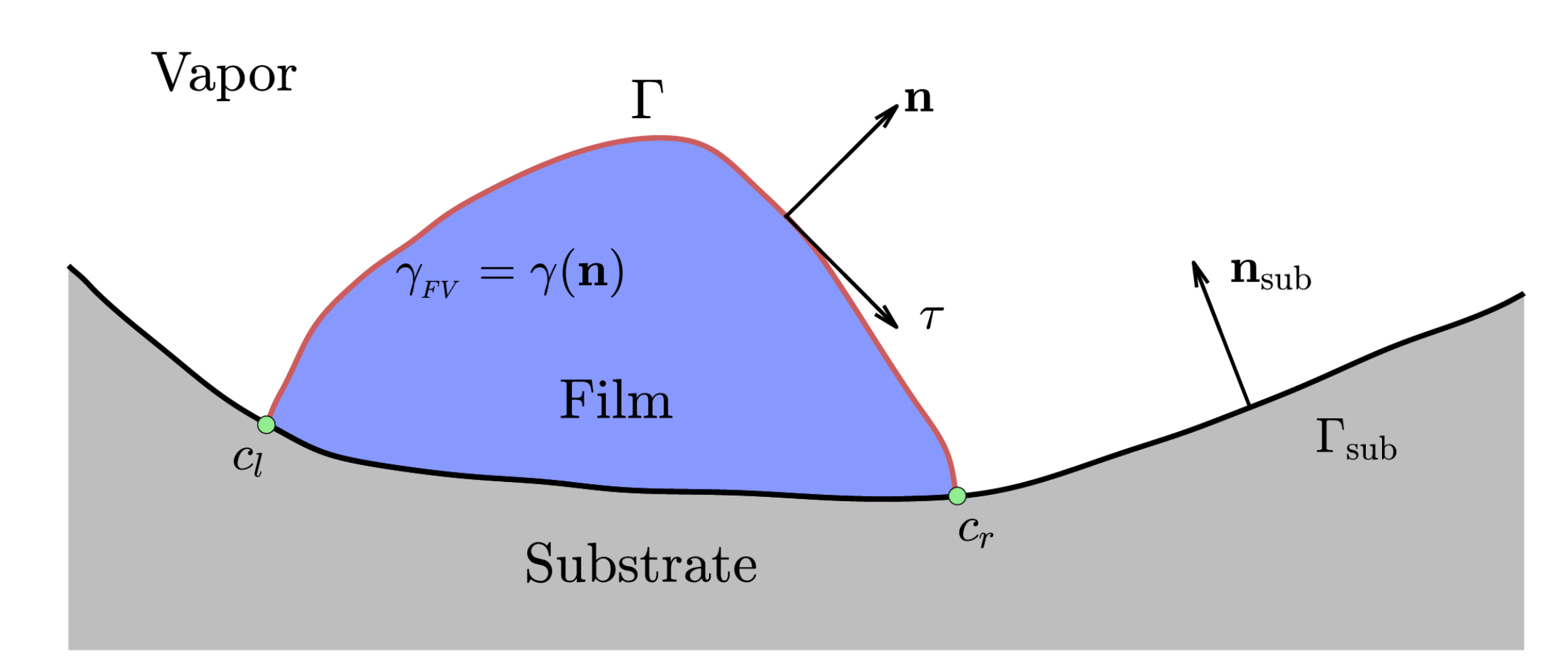}
\caption{A schematic illustration of SSD with anisotropic surface energies on a curved substrate.}
\label{fig:model}
\end{figure}

We consider the 2D case where the film/vapor interface is represented by an evolving curve $\Gamma(t)$ with a parametrization over the reference domain $\bI$ given by
\begin{equation}
\vec{x}(\rho, t)=(x(\rho,t), y(\rho,t))^T: \bar{\bI}\times[0,T]\mapsto\mathbb{R}^2,  \nn
\end{equation}
where $\bI=(0,1)$ is the unit interval with boundary $\partial\bI=\{0,1\}$. We denote by $s$ the arclength of the curve $\Gamma(t)$, which satisfies $\partial_s=|\vec{x}_\rho|^{-1}\partial_\rho$. Throughout the work we assume that $|\vec{x}_\rho(\rho, t)|>0$ in $\bI$. We also introduce the unit tangential and normal vectors
\begin{equation}
\boldsymbol{\tau}(\rho, t) = \vec{x}_s = |\vec{x}_\rho|^{-1}\,\vec{x}_\rho,\qquad \vec{n}(\rho,t) = -(\vec{x}_s)^\perp = -|\vec{x}_\rho|^{-1}(\vec{x}_\rho)^\perp,\nn
\end{equation}
where $(\cdot)^\perp$ represents clockwise rotation of a vector by $\frac{\pi }{2}$.

As shown in Fig.~\ref{fig:model}, the substrate was considered as a fixed smooth curve $\Gamma_{\rm sub}$ which was parameterized via its \textit{arclength parameter} $c$ as
 \begin{equation}
 \vec{r}(c)=(r_1(c), r_2(c))^T: [0, L_{\rm sub}]\mapsto \mathbb{R}^2,\nn
 \end{equation}
where $L_{\rm sub}$ is the total length of $\Gamma_{\rm sub}$. By definition, we know that $|\vec{r}_c| = 1$. And we can thus introduce the unit normal of $\Gamma_{\rm sub}$ as 
\begin{equation}
\vec{n}_{\rm sub}(c) = -(\vec{r}_c)^\perp.\nn
\end{equation}

We further assume that the two endpoints of $\Gamma(t)$ are attached to the curve substrate. This yields two contact points which we denote by $c_l(t), c_r(t)$ in the name of the arclength of $\Gamma_{\rm sub}$. Then the total free energies of the system (scaled by $\gamma_0$)
\begin{equation}\label{eq:energy}
W\left(\vec{x}(\cdot, t)\right) = \int_{\Gamma(t)}\gamma(\vec{n})\,\rd s -\sigma[c_r(t)-c_l(t)], 
\end{equation}
where the material constant $\sigma = \frac{\gamma_{_{\rm VS}}- \gamma_{_{\rm FS}}}{\gamma_0}$ with $\gamma_{_{\rm FS}}$ and $\gamma_{_{\rm VS}}$ representing the film/substrate and vapor/substrate surface energy densities, respectively. Besides, the surface energy density of the film/varpor interface is modeled by an anisotropy function $\gamma(\cdot)$, which belongs to $ C^0(\bR^2)\cap C^2(\bR^2\setminus\{\vec{0}\})$ and satisfies $\gamma>0$ on $\bR^2\setminus\{\vec{0}\}$. We further assume $\gamma(\cdot)$ is positively homogeneous of degree one, 
\begin{equation}\label{eq:homoani}
\gamma(\lambda\vec{p}) = \lambda\gamma(\vec{p})\quad\forall\lambda>0,\quad\vec{p}\in\bR^2.
\end{equation}
This implies that $\gamma(\vec 0) = 0$ and the gradient of $\gamma(\vec p)$ satisfies 
\begin{equation}\label{eq:gammapro}
\nabla\gamma(\vec{p})\cdot\vec{p} = \gamma(\vec{p}),\quad \nabla\gamma(\lambda\vec{p}) = \nabla\gamma(\vec{p}),\qquad\forall\lambda>0, \quad \vec{p}\in\bR^2\setminus\{\vec{0}\} ,
\end{equation}
where $\nabla\gamma$ is the gradient of the anisotropy function $\gamma(\cdot)$.

By the thermodynamic variation of the total energy \eqref{eq:energy}, a dimensionless sharp-interface model for SSD on curved substrates in 2D can be developed (see \cite{Jiang18curved, Jiang19a})
\begin{subequations}\label{eqn:model}
\begin{align}\label{eq:model1}
&\vec{x}_t = \mu_{ss}\,\vec{n},\quad 0< s < L(t),\\
&\mu= -[\nabla\gamma(\vec{n})]_s^\perp\,\cdot\vec{n},\label{eq:model2}
\end{align}
\end{subequations}
where $\mu(\rho,t)$ is the chemical potential or the weighted mean curvature, and the term $\nabla\gamma(\vec{n})$ represents the well-known Cahn-Hoffman vector \cite{Hoffman72, cahn1974vector}. The above equations are supplemented with the following boundary conditions at the two contact points
\begin{itemize}
\item attachment condition
\begin{subequations}\label{eqn:bd}
\begin{equation}\label{eq:bd1}
\vec{x}(0,t) = \vec{r}(c_l),\qquad \vec{x}(1,t) = \vec{r}(c_r),\qquad t\geq 0;
\end{equation}
\item relaxed contact angle condition
\begin{equation}\label{eq:bd2}
\frac{\rd c_l(t)}{\rd t} = \eta\bigl[\nabla\gamma\left(\vec{n}(0,t)\right)\cdot\vec{n}_{\rm sub}(c_l) - \sigma\bigr], \qquad \frac{\rd c_r(t)}{\rd t} = -\eta\bigl[\nabla\gamma\left(\vec{n}(1,t)\right)\cdot\vec{n}_{\rm sub}(c_r) - \sigma\bigr],\qquad t\geq 0;
\end{equation}
\item zero mass-flux condition
\begin{equation}\label{eq:bd3}
\mu_s(0,t) = 0, \qquad \mu_s(1,t) = 0,\qquad t\geq 0.
\end{equation}
\end{subequations}
\end{itemize}
Here \eqref{eq:bd1} ensures that the endpoints of the evolving curve $\Gamma(t)$ are always attached to the substrate curve $\Gamma_{\rm sub}$. Condition in \eqref{eq:bd2} can be interpreted as the relaxed contact angle condition with a prescribed mobility $\eta\in(0,+\infty)$. In the limit $\eta\to +\infty$, they collapse to the anisotropic Young's equation
\begin{equation}
\nabla\gamma(\vec{n})\cdot\vec{n}_{\rm sub} = \sigma\qquad\mbox{at}\quad c_l,c_r,\nn
\end{equation}
which yields a contact angle $\arccos\left(\frac{\sigma}{|\nabla\gamma(\vec{n})|}\right)$ between $\nabla\gamma(\vec{n})$ and $\vec{n}_{\rm sub}$ at the two contact points if $|\nabla\gamma(\vec{n})|\geq \sigma$.  The last condition \eqref{eq:bd3} ensures mass conservation of the thin film. In fact, at the left contact point, differentiating \eqref{eq:bd1} and recalling \eqref{eq:bd2} give rise to 
\begin{equation}
\label{eq:tangentialmotion}
\frac{\rd \vec{x}(0,t)}{\rd t} = \frac{\rd \vec{r}(c_l)}{\rd t} = \vec{r}_c(c_l)\,\frac{\rd c_l(t)}{\rd t} = \eta\bigl[\nabla\gamma\left(\vec{n}(0,t)\right)\cdot\vec{n}_{\rm sub}(c_l) - \sigma\bigr]\,\vec{r}_c(c_l).
\end{equation}
A similar evolving equation can be obtained for the right contact point, which implies that the two contact points should evolve in the tangential direction of $\Gamma_{\rm sub}$, which in turn ensures the attachment condition \eqref{eq:bd1} if initially the two curves are attached.  

The above equations \eqref{eqn:model}, \eqref{eqn:bd} form a complete model for SSD on the curved substrate with anisotropic surface energies. In fact, by taking the time derivative of the energy in \eqref{eq:energy},  it is not difficult to obtain 
\begin{equation}\label{eq:energyim}
\ddt W(\vec{x}(t)) = \int_{\Gamma(t)} [\nabla\gamma(\vec{n})]^\perp\cdot(\vec{x}_t)_s\,\rd s -\sigma\left[\frac{\rd x_r(t)}{\rd t} - \frac{\rd x_l(t)}{\rd t}\right],
\end{equation}
where for simplicity we use $\vec{x}(t) = \vec{x}(\cdot, t)$. Then using \eqref{eq:model1}, \eqref{eq:model2} and \eqref{eq:bd1}, \eqref{eq:bd2}, we obtain the following energy dissipation law
\begin{equation}
\label{eq:energylaw}
\ddt W(\vec{x}(t)) = - \int_{\Gamma(t)}|\mu_s|^2\,\rd s - \frac{1}{\eta}\left[\left(\frac{\rd x_r(t)}{\rd t}\right)^2 + \left(\frac{\rd x_l(t)}{\rd t}\right)^2\right]\leq 0.
\end{equation}
Let $M(\vec{x}(t))$ denote the mass of the thin films, i.e., the area enclosed by the two curves $\Gamma(t)$ and $\Gamma_{\rm sub}$, which can be expressed as
\begin{equation}\label{eq:def of mass}
  M(\vec{x}(t)) = \frac{1}{2}\int_{\Gamma(t)}\vec{x}\cdot\vec{n}\,\rd s - \frac{1}{2} \int_{c_l(t)}^{c_r(t)}\vec{r}\cdot \vec{n}_{\rm sub}\,\rd c,
\end{equation}
where the minus sign is due to that $\vec{n}_{\rm sub}$ is the inward normal to the enclosed region. 
Then recalling the Reynolds transport theorem, and using \eqref{eq:model1} and \eqref{eq:bd3}, we have the mass preservation law
\begin{equation}
\ddt M(\vec{x}(t))=\int_{\Gamma(t)}\vec{x}_t\cdot\vec{n} \rd s - \int_{c_l(t)}^{c_r(t)}\vec{0}\cdot\vec{n}_{\rm sub}\,\rd c= \int_{\Gamma(t)}\mu_{ss}\,\rd s \equiv 0,\quad t\ge0.
\label{eq:masslaw}
\end{equation}

The main aim of this work is to propose a structure-preserving parametric finite element method for the introduced model. This means that we would like to establish an unconditional stability estimate in terms of the energy and an exact mass preservation for the discrete solutions, which mimics \eqref{eq:energylaw} and \eqref{eq:masslaw}, respectively.

\setcounter{equation}{0}
\section{The parametric weak formulation}
\label{sec:weak}

Based on the work in \cite{BL22}, we now introduce a symmetrized surface energy matrix as
\begin{equation}\label{eq:zk}
\mathbf{Z}_k(\vec{n}) = \gamma(\vec{n})I_2 - \vec{n}\,\nabla\gamma(\vec{n})^T - \nabla\gamma(\vec{n})\,\vec{n}^T +k(\vec{n})\,\vec{n}\,\vec{n}^T, \quad \vec{n}\in\mathbb{S}^1,
\end{equation}
where $I_2\in \bR^{2\times 2}$ is the $2\times 2$ identity matrix, and $k(\vec{n})$ is a stabilizing function to ensure that the symmetric matrix $\mathbf{Z}_k(\vec{n})$ satisfies a stability estimate, see Lemma \ref{lem:kZ}.  We have the following lemma for our introduced matrix. 

\begin{lemma}
It holds that 
\begin{subequations}
\begin{align}
\label{eq:bzform}
&\mathbf{Z}_k(\vec{n})\vec{x}_s = [\nabla\gamma(\vec{n})]^\perp,
\\
\label{eq:munform}
&\mu\,\vec{n} = -[\mathbf{Z}_k(\vec{n})\,\vec{x}_s]_s.
\end{align}
\end{subequations}
\end{lemma}
\begin{proof}
On recalling \eqref{eq:gammapro},  we  have a natural decomposition of the Cahn-Hoffman vector $\nabla\gamma(\vec{n})$ as \cite{Hoffman72, Jiang19a} 
\begin{equation}\label{eq:xidecom}
\nabla\gamma(\vec{n}) = (\nabla\gamma(\vec{n})\cdot\vec{n})\,\vec{n} + (\nabla\gamma(\vec{n})\cdot\boldsymbol{\tau})\,\boldsymbol{\tau} = \gamma(\vec{n})\,\vec{n}+ (\nabla\gamma(\vec{n})\cdot\boldsymbol{\tau})\,\boldsymbol{\tau}.\nn
\end{equation}
 Then it is not difficult to show that
\begin{align}
\mathbf{Z}_k(\vec{n})\vec{x}_s &= \mathbf{Z}_k(\vec{n})\,\boldsymbol{\tau} = (\gamma(\vec{n})I_2 - \vec{n}\,\nabla\gamma(\vec{n})^T - \nabla\gamma(\vec{n})\,\vec{n}^T +k(\vec{n})\,\vec{n}\,\vec{n}^T)\,\boldsymbol{\tau}\nn\\
&=  \gamma(\vec{n})\,\boldsymbol{\tau} - (\nabla\gamma(\vec{n})\cdot\boldsymbol{\tau})\,\vec{n} = [\nabla\gamma(\vec{n})]^\perp.\label{eq:xi1}
\end{align}
 On the other hand, we can reformulate \eqref{eq:model2} as \cite{Barrett07Ani, Jiang19a}
\begin{equation}\label{eq:xi2}
\mu\,\vec{n} = -[\nabla\gamma(\vec{n})]_s^\perp.
\end{equation}
Combining \eqref{eq:xi1} and \eqref{eq:xi2} then yields \eqref{eq:munform} directly.
\end{proof}

In fact, in order for the anisotropic surface diffusion in \eqref{eqn:model} to be well-posed, the smooth anisotropy function is required to be convex \cite{Burger07level, Jiang19a}, meaning that 
\begin{equation}\label{eq:convex}
\gamma(\vec u) - \gamma(\vec v) \leq \nabla\gamma(\vec u)\cdot(\vec u-\vec v),\quad\forall\vec u,\vec v\in\bR^2\setminus\{\vec 0\}.
\end{equation}
This helps to a stability proof for a fully implicit scheme, see \cite{Barrett2011}. For arbitrary nozero vectors $\vec h$ and $\hah$, we now set $\vec n = -\frac{\vec h^\perp}{|\vec h|}$ and $\han  = -\frac{\hah^\perp}{|\hah|}$. We then choose $\vec u = -\hah^\perp, \vec v = -\vec h^\perp$ in \eqref{eq:convex}, and using the property of $\nabla\gamma$ in \eqref{eq:gammapro} to obtain that
\begin{equation}
\gamma(-\hah^\perp) - \gamma(-\vec h^\perp)\leq \nabla\gamma(\han)\cdot(-\hah^\perp+\vec h^\perp) =[\nabla\gamma(\han)]^\perp\cdot(\hah-\vec h).\label{eq:gammacc}
\end{equation}
Now using \eqref{eq:bzform} into \eqref{eq:gammacc} we get
\begin{equation}\label{eq:gammacc1}
\gamma(-\hah^\perp) - \gamma(-\vec h^\perp)\leq \left\{\mathbf{Z}_k(\han)\frac{\hah}{|\hah|}\right\}\cdot(\hah-\vec h),
\end{equation}
which holds for arbitrary $k(\vec n)$ in $\bZ_k(\vec n)$ if $\gamma(\cdot)$ is convex.

So far we haven't yet touched upon the role of the stabilising function $k(\vec n)$ in $\bZ_k(\vec n)$. In fact, our next lemma shows that under a very weak assumption on $\gamma(\cdot)$, we can choose a suitable stabilizing function to get an estimate for the introduced matrix $\bZ_k(\vec n)$, which is very similar to the convex condition in \eqref{eq:gammacc1}.
\begin{lemma}\label{lem:kZ}
If the anisotropy function $\gamma(\cdot)$ satisfies 
\begin{equation}\label{eq:zkassump}
\gamma(-\vec n) <3\,\gamma(\vec n),\qquad \forall\vec n\in\mathbb{S}^1,
\end{equation}
then the following inequality holds 
\begin{equation}
\gamma(-\hah^\perp) -\gamma(-\vec h^\perp) \leq  \left\{\mathbf{Z}_k\left(\vec n\right)\frac{\hah}{|\vec h|}\right\}\cdot(\hah - \vec h),\qquad\forall\hah,\vec h\in\bR^2\setminus\{\vec 0\},\label{eq:newconvex}
\end{equation}
for sufficiently large $k(\vec n )$, where $\vec n = -\frac{\vec h^\perp}{|\vec h|}$.
\end{lemma}

\begin{proof}
The proof can be inspired by \cite[Theorem 2.1]{Bao24structure} and here we give a sketch. We first consider two non-zero vectors $\vec{h}, \hat{\vec{h}} \in \mathbb{R}^2$, and denote by $\vec{n} = -\vec{h}^\perp/|\vec{h}|$ and $\hat{\vec{n}} = -\hat{\vec{h}}^\perp/|\hat{\vec{h}}|$. Using vector decomposition and the definition of $\mathbf{Z}_k(\vec{n})$ in \eqref{eq:zk}, we have
  \begin{align}\label{eq:auxiliary sym1}
    \frac{1}{|\vec{h}|}\left(\vec{Z}_{k}(\vec{n}) \hat{\vec{h}}\right)\cdot \hat{\vec{h}}&= \frac{|\hat{\vec{h}}|^2}{|\vec{h}|}\left((\gamma(\vec{n})I_2 -\vec{n}[\nabla\gamma(\vec n)]^T - [\nabla\gamma(\vec n)]\vec{n}^T + k(\vec{n})\vec{n}\vec{n}^T)\hat{\vec{n}}^\perp\right)\cdot \hat{\vec{n}}^\perp \nn\\
    &= \frac{|\hat{\vec{h}}|^2}{|\vec{h}|} \left(\gamma(\vec{n}) - 2(\vec{n}\cdot \hat{\vec{n}}^\perp)(\nabla\gamma(\vec n)\cdot \hat{\vec{n}}^\perp) + k(\vec{n}) (\vec{n}\cdot \hat{\vec{n}})^2\right)\nn\\
    & = \frac{|\hat{\vec{h}}|^2}{|\vec{h}|} \left(-\gamma(\vec{n}) + 2(\vec{n}\cdot \hat{\vec{n}})(\nabla\gamma(\vec n)\cdot \hat{\vec{n}}) + k(\vec{n}) (\vec{n}\cdot \hat{\vec{n}})^2\right) \nn\\
    &:= \frac{|\hat{\vec{h}}|^2}{4|\vec{h}|\gamma(\vec{n})}P_{k}^2(\vec{n}, \hat{\vec{n}}),
  \end{align}
  and
  \begin{align}\label{eq:auxiliary sym2}
    \frac{1}{|\vec{h}|}\left(\vec{Z}_{k}(\vec{n}) \hat{\vec{h}}\right)\cdot \vec{h} &=\frac{|\hat{\vec{h}}||\vec{h}|}{|\vec{h}|}\left((\gamma(\vec{n})I_2 -\vec{n}[\nabla\gamma(\vec n)]^T - [\nabla\gamma(\vec n)]\vec{n}^T + k(\vec{n})\vec{n}\vec{n}^T)\hat{\vec{n}}^\perp\right)\cdot \vec{n}^\perp \nn\\
    &=|\hat{\vec{h}}|\left(\nabla\gamma(\vec n)\cdot \vec{n})(\hat{\vec{n}}\cdot \vec{n}) - (\nabla\gamma(\vec n)\cdot \vec{n}^\perp)(\vec{n}\cdot \hat{\vec{n}}^\perp)\right)\nn\\
    &= |\hat{\vec{h}}|\left((\nabla\gamma(\vec n)\cdot \vec{n})(\hat{\vec{n}}\cdot \vec{n}) + (\nabla\gamma(\vec n)\cdot \vec{n}^\perp)(\hat{\vec{n}}\cdot \vec{n}^\perp)\right)\nn\\
    &=|\hat{\vec{h}}|(\nabla\gamma(\vec n)\cdot \hat{\vec{n}}):=|\hat{\vec{h}}|(Q(\vec{n}, \hat{\vec{n}}) - \gamma(\hat{\vec{n}})).
  \end{align}
  where we introduced the following two auxiliary functions $P_{\alpha}(\vec{n}, \hat{\vec{n}})$ and $Q(\vec{n}, \hat{\vec{n}})$  as
  \begin{subequations}\label{eq:auxiliary sym}
  \begin{align}\label{eq:auxiliary P}
  P_{\alpha}(\vec{n}, \hat{\vec{n}})&:=2\sqrt{(-\gamma(\vec{n})+2(\vec{n}\cdot \hat{\vec{n}})(\nabla\gamma(\vec n)\cdot \hat{\vec{n}})+\alpha (\hat{\vec{n}}\cdot \vec{n}^\perp)^2)\gamma(\vec{n})},\\
  \label{eq:auxiliary Q}
  Q(\vec{n}, \hat{\vec{n}})&:=\gamma(\hat{\vec{n}}) + \nabla\gamma(\vec n)\cdot \hat{\vec{n}}, \qquad\qquad \forall \vec{n}, \hat{\vec{n}}\in\mathbb{S}^1, \, \alpha \geq 0.
  \end{align}
  \end{subequations}

 We next aim to show that $P_{\alpha}(\vec{n}, \hat{\vec{n}})- Q(\vec{n}, \hat{\vec{n}})\geq 0$ for all $\vec{n}, \hat{\vec{n}}\in\mathbb{S}^1$ with sufficiently large $\alpha$. In the case when $\hat{\vec{n}}=-\vec{n}$, by condition \eqref{eq:zkassump}, we have 
  \begin{equation*}
    P_{0}(\vec{n}, -\vec{n}) -Q(\vec{n}, -\vec{n}) = 2\gamma(\vec{n}) - (\gamma(-\vec{n}) - \gamma(\vec{n})) = 3\gamma(\vec{n}) - \gamma(-\vec{n})>0.
  \end{equation*}
 By the continuity of $P_\alpha$ and $Q$, there exists a $0<\delta<1$ such that for all $|\hat{\vec{n}} + \vec{n}|\leq \delta$ we still have $P_{\alpha}(\vec{n}, \hat{\vec{n}})-Q(\vec{n}, \hat{\vec{n}}) \geq P_0(\vec{n}, \hat{\vec{n}})-Q(\vec{n}, \hat{\vec{n}})>0$. For the case $|\hat{\vec{n}} + \vec{n}|>\delta$, we know $|\hat{\vec{n}} - \vec{n}|<2-\delta$. Therefore, in a similar manner to \cite[Theorem 2.1]{Bao24structure}, we get
  \begin{equation}
    (\hat{\vec{n}}\cdot \vec{n}^\perp)^2= \frac{|\vec{n}+\hat{\vec{n}}|^2}{4}|\vec{n}-\hat{\vec{n}}|^2\geq \frac{\delta^2}{4}|\vec{n}-\hat{\vec{n}}|^2, \qquad \gamma(\hat{\vec{n}})\leq \gamma(\vec{n}) + \nabla\gamma(\vec n)\cdot (\hat{\vec{n}} - \vec{n}) + C_1|\hat{\vec{n}}-\vec{n}|^2.
  \end{equation}
We then have the following inequality 
  \begin{align}
    &P_{\alpha}(\vec{n}, \hat{\vec{n}})- Q(\vec{n}, \hat{\vec{n}})= \left(P_{\alpha}(\vec{n}, \hat{\vec{n}})-2(\nabla\gamma(\vec n)\cdot \hat{\vec{n}})\right) - \left(Q(\vec{n}, \hat{\vec{n}})-2(\vec{\xi}\cdot \hat{\vec{n}})\right)\nn\\
    &\geq \left(4\frac{(-\gamma(\vec{n})+2(\vec{n}\cdot \hat{\vec{n}})(\nabla\gamma(\vec n)\cdot \hat{\vec{n}})+\alpha (\hat{\vec{n}}\cdot \vec{n}^\perp)^2)\gamma(\vec{n})- (\nabla\gamma(\vec n)\cdot \hat{\vec{n}})^2}{P(\vec{n}, \hat{\vec{n}}) + 2(\nabla\gamma(\vec n)\cdot \hat{\vec{n}})}\right) \nn\\
    &\qquad\qquad - \left(\nabla\gamma(\vec n)\cdot \hat{\vec{n}} + \gamma(\vec{n}) + \nabla\gamma(\vec n)\cdot (\hat{\vec{n}} - \vec{n}) - 2(\nabla\gamma(\vec n)\cdot \hat{\vec{n}}) + C_1|\vec{n}-\hat{\vec{n}}|^2\right)\nn\\
    &= 4 \frac{-(\nabla\gamma(\vec n)\cdot \vec{n} - \nabla\gamma(\vec n)\cdot \hat{\vec{n}})^2 - \gamma(\vec{n})(\nabla\gamma(\vec n)\cdot \hat{\vec{n}})|\vec{n}-\hat{\vec{n}}|^2 + \alpha \gamma(\vec{n}) (\hat{\vec{n}}\cdot \vec{n}^\perp)^2}{P(\vec{n}, \hat{\vec{n}}) + 2(\nabla\gamma(\vec n)\cdot \hat{\vec{n}})} - C_1|\vec{n}-\hat{\vec{n}}|^2\nn\\
    &\geq \left(4 \frac{-|\nabla\gamma(\vec n)|^2-|\nabla\gamma(\vec n)|\gamma(\vec{n}) + \frac{\delta^2}{4}\alpha \gamma(\vec{n}) }{P(\vec{n}, \hat{\vec{n}}) + 2(\nabla\gamma(\vec n)\cdot \hat{\vec{n}})} - C_1\right)|\vec{n}-\hat{\vec{n}}|^2.
  \end{align}
  Therefore, we can choose $k(\vec{n})$ sufficiently large such that $P_{k(\vec{n})}(\vec{n}, \hat{\vec{n}})- Q(\vec{n}, \hat{\vec{n}})\geq 0$ for all $\vec{n}, \hat{\vec{n}}\in\mathbb{S}^1$. This, together with \eqref{eq:auxiliary sym1} and \eqref{eq:auxiliary sym2}, implies that 
  \begin{align}\label{eq:auxiliary local}
    \frac{1}{|\vec{h}|}\left(\vec{Z}_{k}(\vec{n}) \hat{\vec{h}}\right)\cdot (\hat{\vec{h}}-\vec{h}) &= \frac{|\hat{\vec{h}}|^2}{4|\vec{h}|\gamma(\vec{n})}P_{k(\vec{n})}^2(\vec{n}, \hat{\vec{n}}) - |\hat{\vec{h}}|(Q(\vec{n}, \hat{\vec{n}}) - \gamma(\hat{\vec{n}})) \nn\\
    & \geq |\hat{\vec{h}}| P_{k(\vec{n})}(\vec{n}, \hat{\vec{n}}) - |\hat{\vec{h}}| \gamma(\vec{n}) - |\hat{\vec{h}}| (P_{k(\vec{n})}(\vec{n}, \hat{\vec{n}}) - \gamma(\hat{\vec{n}}))\nn\\
    & = |\hat{\vec{h}}|\gamma(\han) - |\vec{h}| \gamma(\vec{n}) = \gamma(-\hah^\perp) - \gamma(-\vec h^\perp),
  \end{align}
  which gives \eqref{eq:newconvex} as claimed. 
\end{proof}

\begin{rem}
In our fully discrete scheme, we will employ a suitable time splitting of the Cahn-Hoffman vector $\nabla\gamma(\vec n)$ in \eqref{eq:bzform}.  This gives rise to a linear approximation, where an unconditional stability in terms of the anisotropic surface energy will be guaranteed by \eqref{eq:newconvex} naturally, see the proof of Theorem \ref{thm:unstab}. 
\end{rem}
\begin{rem}In fact our assumption on the anisotropy function in \eqref{eq:zkassump} is pretty weak. In most cases $\gamma(\cdot)$ is assumed to be absolutely homogeneous of degree one, i.e., \begin{equation}
\gamma(\vec p) =\gamma(-\vec p)\qquad\forall\vec p\in\bR^2,
\end{equation}
which satisfies \eqref{eq:zkassump} straightforwardly. This case was carefully studied in \cite{BJY21symmetrized} for the evolution of a closed curve by a different mathematical proof framework via global geometric inequality. Based on the
result in this paper, the condition in \cite[(1.13)]{BJY21symmetrized} can be relaxed to \eqref{eq:zkassump} too! It is also noteworthy that we don't require the convexity of $\gamma(\cdot)$ to get \eqref{eq:newconvex}. 
\end{rem}

We introduce the following function space which relies on the two contact point $c_l, c_r$ and $\vec{r}(c)$
\begin{equation}
V_{c_l, c_r} := \left\{\vec{\zeta}\in [H^1(\bI)]^2: \vec{\zeta}(0)\cdot\vec{n}_{\rm sub}(c_l)=0,\; \vec{\zeta}(1)\cdot\vec{n}_{\rm sub}(c_r)=0\right\}.
\end{equation}
Multiplying \eqref{eq:munform} with an arbitrary $\vec{\zeta}\in V_{c_l, c_r}$ and using integration by parts, we get
\begin{align}
0&=\int_{\Gamma(t)}\mu\,\vec{n}\cdot\vec{\zeta}\,\ds - \int_{\Gamma(t)}\left(\mathbf{Z}_k(\vec{n})\,\vec{x}_s\right)\cdot\vec{\zeta}_s\,\ds + \left([\mathbf{Z}_k(\vec{n})\vec{x}_s]\cdot\vec{\zeta}\right)\big|_{\rho=0}^{\rho=1}\nn\\
&= \int_{\Gamma(t)}\mu\,\vec{n}\cdot\vec{\zeta}\,\ds - \int_{\Gamma(t)}\left(\mathbf{Z}_k(\vec{n})\,\vec{x}_s\right)\cdot\vec{\zeta}_s\,\ds + ([\nabla\gamma(\vec{n})]^\perp\cdot\vec{\zeta})\big|_{\rho=0}^{\rho=1}.\label{eq:lrc}
\end{align}
At the contact point $c_l$, we note that $\vec{\zeta}(0)\cdot\vec{n}_{\rm sub}(c_l)=0$. Moreover, the arclength parameterization of $\Gamma_{\rm sub}$ implies that $\{\vec{n}_{\rm sub}(c_l), \vec{r}_c (c_l)\}$ forms an orthonormal basis. By vector decomposition of $[\nabla\gamma(\vec{n}(0,t))]^\perp$, we have
\begin{align}
\vec{\zeta}(0)\cdot[\nabla\gamma(\vec{n}(0,t))]^\perp & = \vec{\zeta}(0)\cdot \left(([\nabla\gamma(\vec{n}(0,t))]^\perp\cdot \vec{r}_c(c_l))\vec{r}_c(c_l)+ ([\nabla\gamma(\vec{n}(0,t))]^\perp\cdot \vec{n}_{\rm sub}(c_l))\vec{n}_{\rm sub}(c_l)\right) \nn\\
& = \vec{\zeta}(0)\cdot\vec{r}_c(c_l)\;  [\nabla\gamma(\vec{n}(0,t))]^\perp\cdot\vec{r}_c(c_l) + \vec{\zeta}(0)\cdot\vec{n}_{\rm sub}(c_l)\;  [\nabla\gamma(\vec{n}(0,t))]^\perp\cdot\vec{n}_{\rm sub}(c_l)\nn\\
& = \vec{\zeta}(0)\cdot\vec{r}_c(c_l)\,[\nabla\gamma(\vec{n}(0,t))]\cdot\vec{n}_{\rm sub}(c_l)\nn\\
&= \vec{\zeta}(0)\cdot\vec{r}_c(c_l)\,\left(\sigma + \frac{1}{\eta}\frac{\rd c_l(t)}{\rd t}\right),\label{eq:leftc}
\end{align}
where we used the condition \eqref{eq:bd2}. In a similar manner, we get for the right contact point that
\begin{equation}
\vec{\zeta}(1)\cdot[\nabla\gamma(\vec{n}(1,t))]^\perp=\vec{\zeta}(1)\cdot\vec{r}_c(c_r)\,\left(\sigma-\frac{1}{\eta}\frac{\rd c_r(t)}{\rd t}\right). 
\label{eq:rightc}
\end{equation}
Using \eqref{eq:leftc} and \eqref{eq:rightc} in \eqref{eq:lrc} gives 
\begin{align}
0&=\int_{\Gamma(t)}\mu\,\vec{n}\cdot\vec{\zeta}\,\ds - \int_{\Gamma(t)}\left(\mathbf{Z}_k(\vec{n})\,\vec{x}_s\right)\cdot\vec{\zeta}_s\,\ds+\sigma\left[\frac{\vec{r}_c(c_r(t))}{|\vec{r}_c(c_r(t))|^2}\cdot\vec{\zeta}(1) - \frac{\vec{r}_c(c_l(t))}{|\vec{r}_c(c_l(t))|^2}\cdot\vec{\zeta}(0)\right]\nn\\
&\hspace{1cm} - \frac{1}{\eta}\,\left(\frac{\rd c_r(t)}{\rd t}\,\frac{\vec{r}_c(c_r(t))}{|\vec{r}_c(c_r(t))|^2}\cdot\vec{\zeta}(1) + \frac{\rd c_l(t)}{\rd t}\,\frac{\vec{r}_c(c_l(t))}{|\vec{r}_c(c_l(t))|^2}\cdot\vec{\zeta}(0)\right).
\end{align}
Here we use the property of the arclength parameterization, i.e., $|\vec{r}_c| = 1$.

We next introduce $\big\langle\cdot,\cdot\big\rangle$ as the $L^2$-inner product over $\bI$, then we are ready to propose the following weak formulation for the introduced model in \S \ref{sec:model}. Given the initial curve $\Gamma(0)$ with the two endpoints $\vec{x}(0,0)=\vec{r}(c_l(0))$ and $\vec{x}(1,0)=\vec{r}(c_r(0))$, for $t\in(0,T]$, we find $\vec{x}(\cdot,t)\in[H^1(\bI)]^2$ with $\vec{x}(0,t)=\vec{r}(c_l(t)), \vec{x}(1,t)=\vec{r}(c_r(t))$, $\vec{x}_t\in V_{c_l(t),c_r(t)}$, and $\mu(\cdot, t)\in H^1(\bI)$ such that
\begin{subequations}\label{eqn:weak}
\begin{align}\label{eq:weak1}
&\big\langle\vec{x}_t\cdot\vec{n},~\chi\,|\vec{x}_\rho|\big\rangle + \big\langle\mu_\rho,~\chi_\rho\,|\vec{x}_\rho|^{-1}\big\rangle=0,\qquad\chi\in H^1(\bI),\\[0.5em]
&\big\langle\mu\,\vec{n},~\vec{\zeta}\,|\vec{x}_\rho|\big\rangle - \big\langle\mathbf{Z}_k(\vec{n})\,\vec{x}_\rho,~\vec{\zeta}_\rho\,|\vec{x}_\rho|^{-1}\big\rangle +\sigma\left[\frac{\vec{r}_c(c_r(t))}{|\vec{r}_c(c_r(t))|^2}\cdot\vec{\zeta}(1) - \frac{\vec{r}_c(c_l(t))}{|\vec{r}_c(c_l(t))|^2}\cdot\vec{\zeta}(0)\right]\nn\\
&\hspace{1cm}- \frac{1}{\eta}\,\left(\frac{\rd c_r(t)}{\rd t}\,\frac{\vec{r}_c(c_r(t))}{|\vec{r}_c(c_r(t))|^2}\cdot\vec{\zeta}(1) + \frac{\rd c_l(t)}{\rd t}\,\frac{\vec{r}_c(c_l(t))}{|\vec{r}_c(c_l(t))|^2}\cdot\vec{\zeta}(0)\right) =0,\qquad\forall\vec{\zeta}\in V_{c_l(t),c_r(t)}.\label{eq:weak2}
\end{align}
\end{subequations}
Now we note that the above formulation is pretty nonlinear as the function space $V_{c_l(t),c_r(t)}$ depends on the solution $\vec{x}(\cdot,t)$ itself. In the case of an isotropic energy where $\mathbf{Z}_k(\vec{n}) = I_2$, combined with a flat substrate such that $\vec{r}_c = (0,1)^T$, \eqref{eqn:weak} collapses to the weak formulation introduced in \cite[(2.11)]{Zhao20}. It is noteworthy that in \cite[(2.11)]{Zhao20}, the flat substrate is already parameterized by its arclength!

We have the following theorem which shows the energy dissipation law and mass conservation law are satisfied within the weak formulation.
\begin{thm}Let $(\vec{x}(\rho,t), \mu(\rho,t))$ be a weak solution of  \eqref{eqn:weak}, then it holds that
\begin{subequations}
\begin{align}\label{eq:weakenergylaw}
&\ddt W(\vec{x}(t)) = -\big\langle\mu_\rho,~\mu_\rho\,|\vec{x}_\rho|^{-1}\big\rangle-\frac{1}{\eta}\left(\left(\frac{\rd c_r(t)}{\rd t}\right)^2+\left(\frac{\rd c_l(t)}{\rd t}\right)^2\right)\leq 0,\\
&\ddt M(\vec{x}(t)) \equiv 0,\qquad t\ge0.
\label{eq:weakmasslaw}
\end{align}
\end{subequations}
\end{thm}

\begin{proof}
Using the Reynolds transport theorem and setting $\chi=1$ in \eqref{eq:weak1}, it is straightforward to obtain
\begin{equation}
\ddt M(\vec{x}(t))=\int_{\Gamma(t)}\vec{x}_t\cdot\vec{n}\,\ds =\big\langle\vec{x}_t\cdot\vec{n},~1\big\rangle = \big\langle \mu_\rho,~1_\rho\,|\vec{x}_\rho|^{-1}\big\rangle=0,\nn
\end{equation}
which implies \eqref{eq:weakmasslaw}.

For the energy stability, we recall \eqref{eq:energyim} and \eqref{eq:bzform} to get
\begin{equation}\label{eq:weakel1}
\ddt W(\vec{x}(t))= \int_{\Gamma(t)}(\mathbf{Z}_k(\vec{n})\,\vec{x}_s)\cdot(\vec{x}_t)_s\,\ds - \sigma\left(\frac{\rd c_r(t)}{\rd t} - \frac{\rd c_l(t)}{\rd t}\right).
\end{equation}
Beside, by \eqref{eq:tangentialmotion}, we have for the left contact point and similarly for the right contact point that it holds
\begin{equation}\label{eq:weakel2}
\frac{\rd \vec{x}(0,t)}{\rd t}\cdot\vec{r}_c(c_l) = \frac{\rd c_l(t)}{\rd t}, \quad \frac{\rd \vec{x}(1,t)}{\rd t}\cdot\vec{r}_c(c_r) = \frac{\rd c_r(t)}{\rd t}.
\end{equation}
Now setting $\chi=\mu$ in \eqref{eq:weak1} and $\vec{\zeta}=\vec{x}_t\in V_{c_l(t),c_r(t)}$ in \eqref{eq:weak2} and using \eqref{eq:weakel1}, \eqref{eq:weakel2}, we obtain \eqref{eq:weakenergylaw} as claimed. 
\end{proof}

\setcounter{equation}{0}
\section{Finite element approximations}\label{sec:fem}

We consider a partition of the time interval as  $[0,T]=\bigcup_{j=0}^{M-1}[t_m, t_{m+1}]$ with $t_m= m\ttau$, where $\ttau = T/M$ is the uniform time step. We further uniformly divide the reference domain as $\bI=\bigcup_{j=1}^{N}\bI_j=\bigcup_{j=1}^N[q_{j-1},q_{j}]$ with $q_j= jh,\, h = 1/N$.  We introduce the finite element spaces
\begin{subequations}
\begin{align}
\mathcal{V}^h&:=\left\{\chi\in C^0(\bI): \chi\big|_{\bI_j}\quad\mbox{is affine}\quad j = 1,\ldots, N  \right\},\\
\mathcal{V}_{c_l, c_r}^{m+1}&:=\left\{\vec{\eta}\in [\mathcal{V}^h]^2:\; \vec{\eta}(0)\cdot(\vec{r}(c_l^{m+1})-\vec{r}(c_l^m))^\perp=0,\;\; \vec{\eta}(1)\cdot(\vec{r}(c_r^{m+1})-\vec{r}(c_r^m))^\perp=0\right\},
\label{eq:deltaXspace}
\end{align}
where $c_{l}^m, c_r^m$ are the approximations of the contact points $c_{l}(t_m), c_r(t_m)$, respectively. Moreover, to enforce the attachment condition \eqref{eq:bd1} for the finite element space, we define the subset of $[\mV^h]^2$ which depends on $c_l,c_r$ as
\begin{equation}
\mathcal{S}^h_{c_l,c_r}: = \left\{\vec{\eta}\in [\mV^h]^2: \vec{\eta}(0)=\vec{r}(c_l), \;\;\vec{\eta}(1)=\vec{r}(c_r)\right\}.
\end{equation}
\end{subequations}
We also introduce the mass-lumped $L^2$-inner product $\langle\cdot,\cdot\rangle^h$ via
\begin{equation}
\left\langle \vec{v}, \vec{w} \right\rangle^h = \tfrac12\,h\sum_{j=1}^J 
\left[(\vec{v}\cdot\vec{w})(q_j^-) + (\vec{v} \cdot \vec{w})(q_{j-1}^+)\right],\label{eq:masslumped}
\end{equation}
for any piecewise continuous functions $\vec{v}, \vec{w}$, 
with possible jumps at the nodes $\{q_j\}_{j=1}^J$, where
$u(q_j^\pm)=\underset{\delta\searrow 0}{\lim}\ u(q_j\pm\delta)$.

We approximate $\Gamma(t_m)$ by a polygonal curve $\Gamma^m$ which has contact with $\Gamma_{\rm sub}$ at the two contact points $c_l^m, c_r^m$, for $m=0,\ldots, M$. Then we introduce $\{\vec{X}^m\}_{0\leq m\leq M}$ with $\vec{X}^m\in\mathcal{S}^h_{c_l^m,c_r^m}$ as the parameterizations of $\{\Gamma^m\}_{0\leq m\leq M}$ such that $\Gamma^m=\vec{X}^m(\bar{\bI})$. The approximations of the unit tangent and normal vectors are given by
\begin{equation}
\boldsymbol{\tau}^m = \vec{X}^m_s = \frac{\vec{X}_\rho^m}{|\vec{X}^m_\rho|},\qquad \vec{n}^m = -(\vec{X}^m_s)^\perp =-\frac{(\vec{X}^m_\rho)^\perp}{|\vec{X}^m_\rho|}.\nn
\end{equation}
We also introduce the time-weighted interface normals \cite{Zhao2021}
\begin{equation}\label{eq:twnormal}
\vec{n}^{m+\frac{1}{2}}= -\frac{1}{2\,|\vec{X}_{\rho}^m|}\left(\vec{X}_{\rho}^m + \vec{X}_{\rho}^{m+1}\right)^\perp.
\end{equation}
 
At the two contact points $c_l, c_r$, we use the chain rule $\vec{r}_c = \frac{\vec{r}_t}{c_t}$, and introduce the following approximations:
\begin{subequations}\label{eqn:rappro}
\begin{align}
\frac{\vec{r}_c(c_l^m)}{|\vec{r}_c(c_l^m)|^2}\approx \vec{G}(c_l^m,c_l^{m+1})= \frac{(c^{m+1}_{l}-c_{l}^m)(\vec{r}(c_l^{m+1})-\vec{r}(c_{l}^m))}{|\vec{r}(c_{l}^{m+1})-\vec{r}(c_l^m)|^2},\\
\frac{\vec{r}_c(c_r^m)}{|\vec{r}_c(c_r^m)|^2}\approx \vec{G}(c_r^m, c_r^{m+1})= \frac{(c^{m+1}_{r}-c_{r}^m)(\vec{r}(c_r^{m+1})-\vec{r}(c_{r}^m))}{|\vec{r}(c_{r}^{m+1})-\vec{r}(c_r^m)|^2}.
\end{align}
\end{subequations}

We then propose an unconditionally stable approximation of the weak formulation \eqref{eqn:weak} as follows. Given the initial curve $\vec{X}^0$ with $\vec X^0(0)=\vec{r}(c_l^0)$ and $\vec{X}^0(1)=\vec{r}(c_r^0)$, for $m\geq 0$, we seek for $c_l^{m+1}$, $c_r^{m+1}$, $\vec{X}^{m+1}\in \mathcal{S}^h_{c_l^{m+1}, c_r^{m+1}}$ with $\delta\vec{X}^{m+1}=\vec{X}^{m+1}-\vec{X}^m\in \mV^{m+1}_{c_l, c_r}$ and $\mu^{m+1}\in \mV^h$ such that 
\begin{subequations}\label{eqn:fd}
\begin{align}\label{eq:fd1}
 &\frac{1}{\ttau}\big\langle(\vec X^{m+1}-\vec{X}^m)\cdot\vec{n}^{m+\frac{1}{2}},~\chi^h\,|\vec{X}^m_\rho|\big\rangle^h +\Big\langle\mu^{m+1}_\rho,~\chi^h_\rho\,|\vec{X}^m_\rho|^{-1}\,\Big\rangle=0,\qquad\forall\chi^h\in \mV^h,\\[0.5em]
&\Big\langle\mu^{m+1}\,\vec{n}^{m+\frac{1}{2}},~\vec{\zeta}^h \,|\vec{X}^m_\rho|\Big\rangle^h - \Big\langle \mathbf{Z}_k(\vec{n}^m)\vec{X}^{m+1}_\rho,~\vec{\zeta}_\rho^h\,|\vec{X}^m_\rho|^{-1}\,\Big\rangle + \sigma\left[\vec{G}(c_r^{m}, c_r^{m+1})\cdot\vec{\zeta}^h(1)- \vec{G}(c_l^m, c_l^{m+1})\cdot\vec{\zeta}^h(0)\right]\nn\\
&\hspace{0.5cm}-\frac{1}{\eta\,\ttau}\left[(c_r^{m+1}-c_r^m)\,\vec{G}(c_r^m, c_r^{m+1})\cdot\vec{\zeta}^h(1) + (c_l^{m+1}-c_l^m)\,\vec{G}(c_l^m, c_l^{m+1})\cdot\vec{\zeta}^h(0)\right]=0,\qquad\forall\vec{\zeta}^h\in \mV^{m+1}_{c_l,c_r},\label{eq:fd2}
\end{align}
\end{subequations}
where $\vec{G}(\cdot,\cdot)$ was introduced in \eqref{eqn:rappro} as an approximation of $\vec{r}_c(c_{l,r}(t))$ in order for the stability of the substrate energy. It is noteworthy that the function space $\mV^m_{c_l,c_r}$ depends on the solutions themselves, meaning that we have used an implicit finite element spaces for both the test function and solution. Therefore, when solving the nonlinear system, special care need to be taken due to the complexity and sensitivity of the system,  as discussed in Section \ref{sec:sol}. 

In the following we present an unconditional stability estimate for the introduced method \eqref{eqn:fd} under a weak assumption on the $\gamma(\cdot)$. 
\begin{thm}\label{thm:unstab}
Assume $\gamma(\vec{n})$ satisfies \eqref{eq:zkassump}, when  $k(\vec n)$ is taken sufficiently large, it holds that
\begin{equation}\label{eq:denergylaw}
  \frac{W(\vec{X}^{m+1})-W(\vec{X}^m)}{\ttau}+ \big\langle|\mu_\rho^{m+1}|^2,~|\vec{X}_\rho^m|\big\rangle + \frac{1}{\eta}\left[\left(\frac{c_r^{m+1}-c_r^m}{\ttau}\right)^2+\left(\frac{c_l^{m+1}-c_l^m}
  {\ttau}\right)^2\right]\leq 0,
  \end{equation}
  which is the discrete analogue of the energy law in \eqref{eq:energylaw}, and this immediately implies the
  energy dissipation in the fully discrete level, i.e.
  \begin{equation}
  W(\vec{X}^{m+1})\le W(\vec{X}^m)\le \cdots \le W(\vec{X}^0), \qquad m\ge0.
  \end{equation} 
\end{thm}

 \begin{proof}
  On recalling Lemma \ref{lem:kZ}, we have \eqref{eq:newconvex} for sufficiently large $k(\vec n)$. By taking $\vec{h} = \vec{X}^{m}_\rho$ and $\hat{\vec{h}} = \vec{X}^{m+1}_\rho$ in \eqref{eq:newconvex} on each element $\bI_j$, and then summing over all elements, we have
  \begin{align}
    \Big\langle\mathbf{Z}_k(\vec{n}^m)\vec{X}_\rho^{m+1},~(\vec{X}_\rho^{m+1}-\vec{X}_\rho^m)\,|\vec{X}_\rho^m|^{-1}\Big\rangle \geq \Big\langle\gamma(\vec{n}^{m+1}),~|\vec{X}_\rho^{m+1}|\Big\rangle- \Big\langle\gamma(\vec{n}^m),~|\vec{X}^m_\rho|\Big\rangle.\label{eq:zkde}
  \end{align}

  On the other hand, we choose $\chi=\ttau\,\mu^{m+1}$ in \eqref{eq:fd1} and $\vec{\zeta}^h = \vec{X}^{m+1}-\vec{X}^m \in \mathcal{V}_{c_l, c_r}^{m+1}$ in \eqref{eq:fd2}. Combining these two equations then yields 
 \begin{align}
 &\ttau\,\Big\langle\mu_\rho^{m+1},~\mu_\rho^{m+1}\,|\vec{X}^m_\rho|\,\Big\rangle + \Big\langle\mathbf{Z}_k(\vec{n}^m)\vec{X}_\rho^{m+1},~(\vec{X}_\rho^{m+1}-\vec{X}_\rho^m)\,|\vec{X}_\rho^m|^{-1}\Big\rangle-\sigma[(c_r^{m+1}-c_r^m) - (c_l^{m+1}-c_l^m)]\nn\\
 &\hspace{2cm}+\frac{1}{\eta\,\ttau}[(c_r^{m+1}-c_r^m)^2 + (c_l^{m+1}-c_l^m)^2]=0,
 \label{eq:denergy1}
 \end{align}
 where we used the fact $\vec{X}^{m+1}(0)-\vec{X}^m(0) = \vec{r}(c_l^{m+1})-\vec{r}(c_l^m)$, $\vec{X}^{m+1}(1)-\vec{X}^m(1) = \vec{r}(c_r^{m+1})-\vec{r}(c_r^m)$.

We then insert \eqref{eq:zkde} into \eqref{eq:denergy1} and recall \eqref{eq:energy} to get that 
\begin{equation}
W(\vec{X}^{m+1})-W(\vec{X}^m) +\ttau\,\Big\langle\mu_\rho^{m+1},~\mu_\rho^{m+1}\,|\vec{X}^m_\rho|\,\Big\rangle +\frac{1}{\eta\,\ttau}[(c_r^{m+1}-c_r^m)^2 + (c_l^{m+1}-c_l^m)^2]\leq 0,\nn
\end{equation}
which can be recast as \eqref{eq:denergylaw}.
 \end{proof}
 
  We have the following error estimate for the enclosed mass between the polygonal curve $\Gamma^m$ and the substrate curve $\Gamma_{\rm sub}$, which is an approximation of the mass preservation law in \eqref{eq:masslaw}.
 
\begin{figure}
\centering
\includegraphics[width=0.6\textwidth]{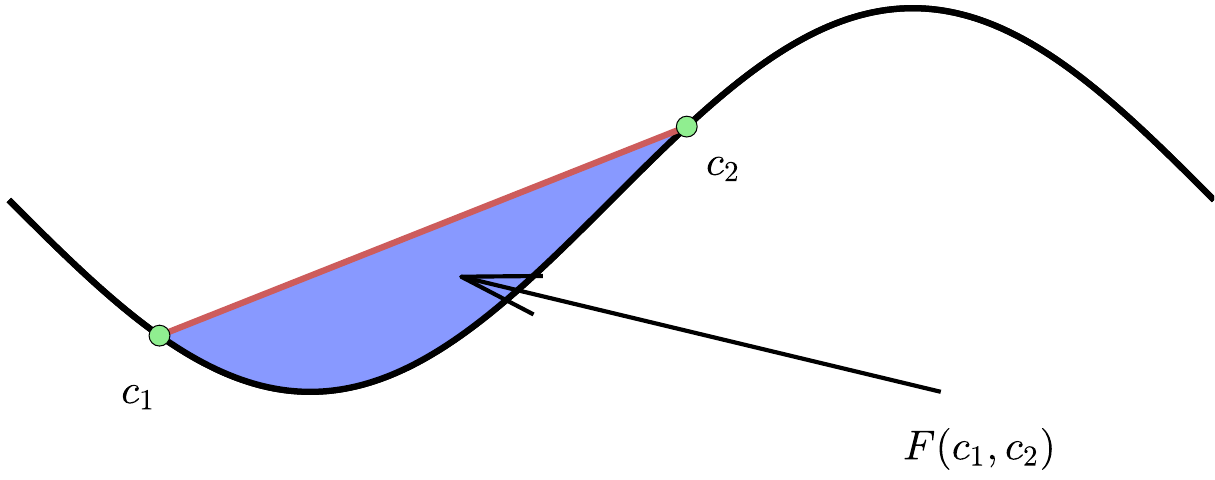}
\caption{Illustration of the area $F(c_1, c_2)$ in \eqref{eq:area}.}
\label{fig:illuF}
\end{figure}

\begin{thm}[mass preservation estimate]
Let $(\vec{X}^{m+1}, \mu^{m+1}, c_l^{m+1}, c_r^{m+1})$ be a solution of the method \eqref{eqn:fd}. Then it holds that
\begin{equation}\label{eq:massperror}
M(\vec{X}^{m+1})-M (\vec{X}^m)= -F(c_l^{m}, c_l^{m+1}) + F(c_r^{m}, c_r^{m+1}),
\end{equation}
where we introduced 
\begin{equation}\label{eq:area}
F(c_1, c_2):=\frac{1}{2}(r_1(c_2)-r_1(c_1))(r_2(c_2)+r_2(c_1))-\frac{1}{2}\int_{c_1}^{c_2} \vec{r}\cdot \vec{n}_{\rm sub}\,\rd c
\end{equation}
as the area between $\Gamma_{\rm sub} $ and the line segment connecting the two points $\vec{r}(c_1),\; \vec{r}(c_2)$, as shown in Fig. \ref{fig:illuF}.  
\end{thm}

 \begin{proof}
  We consider appropriate extensions of $\vec{X}^m, \vec{X}^{m+1}$ to $[0, 4]$ as follows:
  \begin{equation}\label{eq:extension of X}
    \vec{X}^m_{\text{ext}}(\rho) = 
    \begin{cases}
    \vec{X}^m(\rho), & \rho \in [0,1], \\
    \vec{X}^m(1), & \rho \in (1,2], \\
    \vec{X}^m(3-\rho), & \rho \in (2,3], \\
    \vec{X}^m(0), & \rho \in (3,4],
    \end{cases}\qquad \quad 
    \vec{X}^{m+1}_{\text{ext}}(\rho) =
    \begin{cases}
    \vec{X}^{m+1}(\rho), & \rho \in [0,1], \\
    (2-\rho)\vec{X}^{m+1}(1) + (\rho-1)\vec{X}^m(1), & \rho \in (1,2], \\
    \vec{X}^{m}(3-\rho), & \rho \in (2,3], \\
    (4-\rho)\vec{X}^{m}(0) + (\rho-3) \vec{X}^{m+1}(0), & \rho \in (3,4].
    \end{cases}
    \end{equation}
Let $\Gamma^m_{\text{ext}} = \vec{X}^m_{\text{ext}}([0, 4])$ and $\Gamma^{m+1}_{\text{ext}} = \vec{X}^{m+1}_{\text{ext}}([0, 4])$. It is easy to see that $\Gamma^m_{\text{ext}}$ and $\Gamma^{m+1}_{\text{ext}}$ are two closed polygonal curves, and $\Gamma^m_{\text{ext}}$ is degenerated. We denote by $\mathcal{A}^m_{\text{ext}}$ and $\mathcal{A}^{m+1}_{\text{ext}}$ the areas enclosed by $\Gamma^m_{\text{ext}}$ and $\Gamma^{m+1}_{\text{ext}}$, respectively. 

In a similar manner to \cite{Zhao2021}, we introduce the intermediate closed curve between $\Gamma^m_{\text{ext}}$ and $\Gamma^{m+1}_{\text{ext}}$ as
\begin{equation}
\vec{X}_{\text{ext}}^h(\rho, t) := \frac{t_{m+1}-t}{\ttau}\vec{X}^m_{\text{ext}}(\rho) + \frac{t-t_m}{\ttau}\vec{X}^{m+1}_{\text{ext}}(\rho)\quad\mbox{for}\quad  t\in[t_m,t_{m+1}],\;\forall \rho \in [0, 4],
\end{equation}
and denote by $\Gamma^h_{\text{ext}}(t) = \vec{X}_{\text{ext}}^h([0, 4])$. Let $\mathcal{A}_{\text{ext}}(t)$ be the area of $\Gamma^h_{\text{ext}}(t)$. By Reynolds transport theorem, we obtain
\begin{align}\label{eq:ve1}
\ddt\mathcal{A}_{\text{ext}}(t) &= \int_{\Gamma^h_{\text{ext}}(t)} (\vec{X}_{\text{ext}}^h)_t\cdot\vec{n}_{\text{ext}}^h\,\rd s=\frac{1}{\ttau}\int_0^4(\vec{X}_{\text{ext}}^{m+1}-\vec{X}_{\text{ext}}^m)\cdot[-(\vec{X}_{\text{ext}}^h)_{\rho}]^\perp\,\rd\rho\nonumber\\
&= \frac{1}{\ttau}\int_0^1 (\vec{X}^{m+1}-\vec{X}^m)\cdot(-(\vec{X}_{\text{ext}}^h)_{\rho})^\perp\,\rd\rho + \frac{1}{\ttau}\int_1^2 (2-\rho)(\vec{X}^{m+1}-\vec{X}^m)(1)\cdot \left(\frac{t-t_m}{\ttau}(\vec{X}^{m+1}-\vec{X}^{m})(1)\right)^\perp\,\rd\rho\nonumber\\
&\qquad + \frac{1}{\ttau}\int_3^4 (\rho-3)(\vec{X}^{m+1}-\vec{X}^m)(0)\cdot\left(\frac{t-t_m}{\ttau}(\vec{X}^{m}-\vec{X}^{m+1})(0)\right)^\perp\,\rd\rho\nonumber\\
&= \frac{1}{\ttau}\int_0^1 (\vec{X}^{m+1}-\vec{X}^m)\cdot\left(-\frac{t_{m+1}-t}{\ttau}\vec{X}^m_{\rho}(\rho) -\frac{t-t_m}{\ttau}\vec{X}^{m+1}_{\rho}(\rho)\right)^\perp\,\rd\rho,
\end{align}
where we use that fact $\vec{n}^h_{\text{ext}} = -\frac{((\vec{X}_{\text{ext}}^h)_{\rho})^\perp}{|(\vec{X}_{\text{ext}}^h)_{\rho}|}$ and $(\vec{X}_{\text{ext}}^h)_t = \frac{\vec{X}^{m+1}_{\text{ext}}-\vec{X}^m_{\text{ext}}}{\ttau}$. Integrating \eqref{eq:ve1} for $t$ from $t_m$ to $t_{m+1}$ then gives us
\begin{align}
  \mathcal{A}_{\text{ext}}^{m+1}-\mathcal{A}_{\text{ext}}^m &= \int_{t_m}^{t_{m+1}}\frac{1}{\ttau}\int_0^1 (\vec{X}^{m+1}-\vec{X}^m)\cdot\left(-\frac{t_{m+1}-t}{\ttau}\vec{X}^m_{\rho}(\rho) -\frac{t-t_m}{\ttau}\vec{X}^{m+1}_{\rho}(\rho)\right)^\perp\,\rd\rho\, \rd t\nn\\
& = \int_{0}^1(\vec{X}^{m+1}-\vec{X}^m)\cdot\frac{1}{\ttau}\int_{t_m}^{t_{m+1}}\left(-\frac{t_{m+1}-t}{\ttau}\vec{X}^m_{\rho}(\rho) -\frac{t-t_m}{\ttau}\vec{X}^{m+1}_{\rho}(\rho)\right)^\perp\rd t\rd\rho\nn\\
&=\int_{0}^1(\vec{X}^{m+1}-\vec{X}^m)\cdot\vec{n}^{m+\frac{1}{2}}\,\rd \rho = \big\langle (\vec{X}^{m+1}-\vec{X}^m)\cdot\vec{n}^{m+\frac{1}{2}},~1\big\rangle.\label{eq:ve2}
\end{align}

On the other hand, we know $\mathcal{A}_{\text{ext}}^m = 0$, and $\mathcal{A}_{\text{ext}}^{m+1}$ can be given as
\begin{align}
  \mathcal{A}_{\text{ext}}^{m+1} = \frac{1}{2} \int_{\Gamma^{m+1}_{\text{ext}}} \vec{X}_{\text{ext}}^{m+1}\cdot\vec{n}_{\text{ext}}^{m+1}\,\rd s &= \frac{1}{2}\int_{\Gamma^{m+1}} \vec{X}^{m+1}\cdot\vec{n}^{m+1}\,\rd s - \frac{1}{2}\int_{\Gamma^m} \vec{X}^m\cdot\vec{n}^m\,\rd s \nn\\
  &\quad -\frac{1}{2}(x^{m+1}(\rho_N) - x^m(\rho_N))(y^{m+1}(\rho_N)+y^m(\rho_N)) \nn\\
  &\quad + \frac{1}{2}(x^{m+1}(\rho_0) - x^m(\rho_0))(y^{m+1}(\rho_0)+y^m(\rho_0)).
\end{align}
This, together with \eqref{eq:ve2}, the definition of $M$ in \eqref{eq:def of mass}, and the definition of $F$, gives
\begin{align}\label{eq:ve3}
M(\vec{X}^{m+1})-M (\vec{X}^m)&= \frac{1}{2}\int_{\Gamma^{m+1}} \vec{X}^{m+1}\cdot\vec{n}^{m+1}\,\rd s - \frac{1}{2}\int_{\Gamma^m} \vec{X}^m\cdot\vec{n}^m\,\rd s - \frac{1}{2}\int_{c_l^{m+1}}^{c_r^{m+1}} \vec{r}\cdot \vec{n}_{\rm sub}\,\rd c\nn + \frac{1}{2}\int_{c_l^{m}}^{c_r^{m}} \vec{r}\cdot \vec{n}_{\rm sub}\,\rd c\nn\\
& = \big\langle (\vec{X}^{m+1}-\vec{X}^m)\cdot\vec{n}^{m+\frac{1}{2}},~1\big\rangle - F(c_l^{m}, c_l^{m+1}) + F(c_r^{m}, c_r^{m+1}).
\end{align}

Now choosing $\chi^h = 1$ in \eqref{eq:fd1} and combining \eqref{eq:ve3} then yields \eqref{eq:massperror}.
 \end{proof}
 
 \begin{rem}
By the isoperimetric inequality, we note that 
\begin{align}
|F(c_l^{m}, c_l^{m+1})| \leq \frac{1}{4\pi}\left(|\vec{r}(c_l^{m+1})-\vec{r}(c_l^m)| + |c_l^{m+1}-c_l^m|\right)^2 \leq \frac{1}{\pi}|c_l^{m+1} - c_l^m|^2.\label{eq:Fb} 
\end{align}
We thus have 
\begin{equation}
  |F(c_l^{m}, c_l^{m+1})|\lesssim |c_l^{m+1} - c_l^m|^2,\nn
\end{equation}
which combining \eqref{eq:massperror} yields that  
\[|M(\vec{X}^{m+1})- M (\vec{X}^m)|\leq |F(c_l^{m}, c_l^{m+1})| +|F(c_r^{m}, c_r^{m+1})|\lesssim(|c_l^{m+1}-c_l^m|^2 + |c_r^{m+1}-c_r^m|^2).\]
In the case of a flat substrate,  both $F(c_l^{m}, c_l^{m+1}) = F(c_r^{m}, c_r^{m+1}) = 0$, and thus we have an exact mass preservation 
\[M (\vec{X}^{m+1})-M(\vec{X}^m)=0.\]
 \end{rem}

Although we have employed the time-weighted discrete normals in \eqref{eqn:fd}, it is still not possible to enable an exact mass preservation for the area enclosed by $\Gamma^m$ and the curve substrate $\Gamma_{\rm sub}$. We next aim to make a correction to the discrete normal $\vec{n}^{m+\frac{1}{2}}$ so that an exact mass preservation can be achieved. We introduce $\delta \vec{n}^{m+\frac{1}{2}}\in V^h$ such that
\begin{equation}\label{eq:correctnormal}
\delta \vec{n}^{m+\frac{1}{2}}(q_j) =\left\{\begin{array}{ll}
-\frac{2\,F(c_l^{m},c_l^{m+1}) }{|\vec{r}(c_l^{m+1})-\vec{r}(c_l^m)|^2}\frac{(\vec{r}(c_l^{m+1}) - \vec{r}(c_l^m))}{|\vec{X}^m(\rho_1) - \vec{X}^m(\rho_0)|}\quad &\mbox{for}\quad j = 0,\\
\frac{2\,F(c_r^{m},c_r^{m+1}) }{|\vec{r}(c_r^{m+1})-\vec{r}(c_r^m)|^2}\frac{(\vec{r}(c_r^{m+1}) - \vec{r}(c_r^m))}{|\vec{X}^m(\rho_N)-\vec{X}^m(\rho_{N-1})|}\quad &\mbox{for}\quad j = N,\\
0\quad &\mbox{otherwise}.
\end{array}\right.
\end{equation}
 
 We then adapt \eqref{eqn:fd} to obtain the mass preservation property as follows. Given the initial curve $\vec{X}^0$ with $\vec{X}^0(0)=\vec{r}(c_l^0)$ and $\vec{X}^0(1)=\vec{r}(c_r^0)$, for $m\geq 0$, we seek for $c_l^{m+1}$, $c_r^{m+1}$, $\vec{X}^{m+1}\in \mathcal{S}^h_{c_l^{m+1}, c_r^{m+1}}$ with $\delta\vec{X}^{m+1}=\vec{X}^{m+1}-\vec{X}^m\in \mV^{m+1}_{c_l, c_r}$ and $\mu^{m+1}\in \mV^h$ such that
\begin{subequations}\label{eqn:vfd}
\begin{align}\label{eq:vfd1}
&\frac{1}{\ttau}\Big\langle(\vec{X}^{m+1}-\vec{X}^m)\cdot\vec{n}^{m+\frac{1}{2}}_*,~\chi^h\,|\vec{X}^m_\rho|\Big\rangle^h +\Big\langle\mu^{m+1}_\rho,~\chi^h_\rho\,|\vec{X}^m_\rho|^{-1}\,\Big\rangle=0,\qquad\forall\chi^h\in \mV^h,\\[0.7em]
&\Big\langle\mu^{m+1}\,\vec{n}^{m+\frac{1}{2}}_*,~\vec{\zeta}^h \,|\vec{X}^m_\rho|\Big\rangle^h - \Big\langle\mathbf{Z}_k(\vec{n}^m)\vec{X}^{m+1}_\rho,~\vec{\zeta}_\rho^h\,|\vec{X}^m_\rho|^{-1}\,\Big\rangle + \sigma\left[\vec{G}(c_r^{m}, c_r^{m+1})\cdot\vec{\zeta}^h(1)- \vec{{G}}(c_l^m, c_l^{m+1})\cdot\vec{\zeta}^h(0)\right]\nn\\
&\hspace{0.5cm}-\frac{1}{\eta\,\ttau}\left[(c_l^{m+1}-c_l^m)\,\vec{G}(c_l^m, c_l^{m+1})\cdot\vec{\zeta}^h(1) + (c_r^{m+1}-c_r^m)\,\vec{G}(c_r^m, c_r^{m+1})\cdot\vec{\zeta}^h(0)\right]=0,\qquad\forall\vec{\zeta}^h\in \mV^{m+1}_{c_l,c_r},\label{eq:vfd2}
\end{align}
\end{subequations}
where $\vec{n}^{m+\frac{1}{2}}_* := \vec{n}^{m+\frac{1}{2}}+\delta\vec{n}^{m+\frac{1}{2}}$ is the corrected time-weighted interface normals according to the profile of the substrate $\Gamma_{\rm sub}$. 

In the limit of $\ttau\to 0$, we note that from \eqref{eq:correctnormal} that 
\begin{equation}
|\delta \vec{n}^{m+\frac{1}{2}}(q_0)|=\frac{2\,|F(c_l^{m},c_l^{m+1})| }{|\vec{r}(c_l^{m+1})-\vec{r}(c_l^m)|^2}\frac{|\vec{r}(c_l^{m+1}) - \vec{r}(c_l^m)|}{|\vec{X}^m(\rho_1) - \vec{X}^m(\rho_0)|}\lesssim \frac{|\vec{r}(c_l^{m+1}) - \vec{r}(c_l^m)|}{|\vec{X}^m(\rho_1) - \vec{X}^m(\rho_0)|}=O(\ttau/h),
\end{equation}
where we recall \eqref{eq:Fb}. Since $\delta \vec{n}^{m+\frac{1}{2}}$ is corrected locally at $q_0$ and $q_N$, we then note that the mass correction scheme \eqref{eqn:vfd} is a $O(\ttau)$ correction to the unconditional stable scheme in \eqref{eqn:fd}.

We have the following theorem for the new modified method, which mimics both the energy dissipation law \eqref{eq:energylaw} and the mass preservation law \eqref{eq:masslaw} on discrete level. 

\begin{thm}[structure-preserving property]\label{thm:sp}

Let $(\vec{X}^{m+1}, \mu^{m+1}, c_l^{m+1}, c_r^{m+1})$ be a solution of the adapted method \eqref{eqn:vfd}. Assume that the matrix $\mathbf{Z}_k(\vec{n})$ in \eqref{eq:zk} satisfies \eqref{eq:zkassump}, then we have the unconditional energy stability estimate \eqref{eq:denergylaw} for $m=0,\ldots, M-1$. Moreover, it holds for $m=0,\ldots, M-1$
\begin{equation}\label{eq:vdmasslaw}
M(\vec{X}^{m+1})- M(\vec{X}^m)= 0.
\end{equation}
\end{thm}
\begin{proof}
Setting $\chi^h = \ttau\,\mu^{m+1}$ in \eqref{eq:vfd1} and $\vec{\zeta}^h = \vec{X}^{m+1}-\vec{X}^m$ in \eqref{eq:vfd2}, combining these equations yields \eqref{eq:denergy1}. The remain proof is exactly the same as that of the Theorem  \ref{thm:unstab}.

For the mass preservation, we set $\chi^h = \ttau$ in \eqref{eq:vfd1} and obtain that
\begin{equation}
\big\langle(\vec{X}^{m+1} - \vec{X}^m)\cdot\vec{n}^{m+\frac{1}{2}},~1\big\rangle^h = -\big\langle(\vec{X}^{m+1} - \vec{X}^m)\cdot\delta\vec{n}^{m+\frac{1}{2}},~1\big\rangle^h = F( c_l^{m}, c_l^{m+1}) - F( c_r^{m}, c_r^{m+1}) .
\end{equation}
By \eqref{eq:ve2}, \eqref{eq:ve3}, we then obtain \eqref{eq:vdmasslaw} as claimed.
\end{proof}

  \setcounter{equation}{0}
\section{Numerical results}\label{sec:num}

\subsection{Solution method}\label{sec:sol}

For the introduced nonlinear method in \eqref{eqn:vfd}, we introduce an iterative algorithm as the solution method. At each time $t_m$, we compute the solution as follows. Initially, we set $\mu^{m+1,0}=\mu^m$,  $\vec{X}^{m+1, 0}(\rho_j) = \vec{X}^{m}(\rho_j)$ for $j \in {1, \ldots, N-1}$, while on the boundary we set $\vec{X}^{m+1, 0}(\rho_0)=\vec{r}(c_{l}^{m+1, 0})$ and $\vec{X}^{m+1, 0}(\rho_N)=\vec{r}(c_r^{m+1, 0})$. The contact points are initialized with a small random perturbation: $c_{l,r}^{m+1,0} = c_{l,r}^m + 10^{-8}\mathcal{N}(0, 1)$, where $\mathcal{N}(0,1)$ represents the normal distribution. This serves to prevent numerical degeneracy. Then for each $\ell\geq 0$, we find the two contact points $c_l^{m+1, \ell+1}, c_r^{m+1, \ell+1}$, $\vec{X}^{m+1,\ell+1} = \vec{X}^{m}+\delta\vec{X}^{m+1,\ell+1}$ with $\delta\vec{X}^{m+1, \ell+1}\in [\mV^h]^2$ and $\mu^{m+1,\ell+1}\in\mV^h$ such that
\begin{subequations}\label{eqn:vpc}
\begin{align}\label{eq:vpc1}
&\frac{1}{\ttau}\Big\langle(\vec{X}^{m+1, \ell+1}-\vec{X}^m)\cdot\vec{n}^{m+\frac{1}{2}, \ell}_*,~\chi^h\,|\vec{X}^m_\rho|\Big\rangle^h +\Big\langle\mu^{m+1,\ell+1}_\rho,~\chi^h_\rho\,|\vec{X}^m_\rho|^{-1}\,\Big\rangle=0\qquad\forall\chi^h\in \mV^h,\\[0.7em]
&\Big\langle\mu^{m+1,\ell+1}\,\vec{n}^{m+\frac{1}{2},\ell}_*,~\vec{\zeta}^h \,|\vec{X}^m_\rho|\Big\rangle^h - \Big\langle\mathbf{Z}_k(\vec{n}^m)\vec{X}^{m+1, \ell+1}_\rho,~\vec{\zeta}_\rho^h\,|\vec{X}^m_\rho|^{-1}\,\Big\rangle + \sigma\left[\vec{G}(c_r^{m}, c_r^{m+1,\ell})\cdot\vec{\zeta}^h(1)- \vec{G}(c_l^m, c_l^{m+1,\ell})\cdot\vec{\zeta}^h(0)\right]\nn\\
&\hspace{0cm}-\frac{1}{\eta\,\ttau}\left[(c_r^{m+1,\ell+1}-c_r^m)\,\vec{G}(c_r^m, c_r^{m+1,\ell})\cdot\vec{\zeta}^h(1)+(c_l^{m+1,\ell+1}-c_l^m)\,\vec{G}(c_l^m, c_l^{m+1,\ell})\cdot\vec{\zeta}^h(0)\right]=0\quad\forall\vec{\zeta}^h\in \mV^{m+1, \ell}_{c_l,c_r},\label{eq:vpc2}
\end{align}
\end{subequations}
with the linear constraints 
\begin{subequations}\label{eq:vpc3}
\begin{align}
&\delta\vec{X}^{m+1,\ell+1}(\rho_0) = \vec{r}(c_{l}^{m+1, \ell}) + \vec{r}_c(c_{l}^{m+1, \ell})(c_{l}^{m+1, \ell+1}-c_{l}^{m+1, \ell}),\\
&\delta\vec{X}^{m+1,\ell+1}(\rho_N) = \vec{r}(c_{r}^{m+1, \ell}) + \vec{r}_c(c_{r}^{m+1, \ell})(c_{r}^{m+1, \ell+1}-c_{r}^{m+1, \ell}),
\end{align}
\end{subequations}
where $\mV^{m+1, \ell}_{c_l,c_r}$ is defined through \eqref{eq:deltaXspace} except that $c_l^{m+1}, c_r^{m+1}$ are replaced with $c_l^{m+1,\ell}, c_r^{m+1,\ell}$.  Besides, we introduce $\vec{n}^{m+\frac{1}{2},\ell}_*$ as 
\[\vec{n}^{m+\frac{1}{2},\ell}_*=\vec{n}^{m+\frac{1}{2},\ell} + \delta\vec{n}^{m+\frac{1}{2},\ell},\]
where $\vec{n}^{m+\frac{1}{2},\ell}$ and $\delta\vec{n}^{m+\frac{1}{2},\ell}$ follow \eqref{eq:twnormal} and \eqref{eq:correctnormal}, respectively, except that $\vec{X}^{m+1}$ and $c_{l,r}^{m+1}$ are replaced with $\vec{X}^{m+1,\ell}$ and $c_{l,r}^{m+1,\ell}$, respectively. 

In our iterative approach,  the strong restriction $\vec{X}^{m+1} = \vec{X}^m + \delta \vec{X}^{m+1}\in \mathcal{S}^h_{c_l^{m+1}, c_r^{m+1}}$ in the numerical scheme \eqref{eqn:vfd} is relaxed by $\delta \vec{X}^{m+1, \ell+1}\in [\mV^h]^2$ with a Newton type approximation \eqref{eq:vpc3}. We repeat the above iteration until the following condition is satisfied
\begin{equation}
\max_{0\leq j \leq N}|\vec{X}^{m+1,\ell+1}(\rho_j)-\vec{X}^{m+1,\ell}(\rho_j)|\leq {\rm tol},\label{eq:tol}
\end{equation}
where ${\rm tol}$ is a chosen tolerance. Once we have the convergent solution to get $\vec{X}^{m+1}\in\mV^h$, we then simply set $\vec{X}^{m+1}\in\mathcal{S}^h_{c_l^{m+1}, c_r^{m+1}}$ by choosing $\vec{X}^{m+1}(\rho_{0,N})=\vec{r}(c_{l,r}^{m+1})$. 

\begin{rem}
The introduced iterative algorithm in \eqref{eqn:vpc}  and \eqref{eq:vpc3} can be interpreted as a hybrid method, where we employed a Picard-type iteration for the computation of $c_{l,r}^{m+1,\ell+1}$ and $\vec{X}^{m+1, \ell+1}$ at the interior points. While for $\vec{X}^{m+1,\ell+1}$ at the boundary, we approximate the attachment constraint with a Newton's iterative method, and the strong constraint $\vec{X}^{m+1,\ell+1}\in\mathcal{S}^h_{c_l^{m+1}, c_r^{m+1}}$ is only enforced for the final convergent solution. An alternative approach is to update $\vec{X}^{m+1,\ell+1}(\rho_{0,N})=\vec{r}(c_{l,r}^{m+1,\ell+1})$ at each iteration, which however might fail to converge. This is possible due to that the strong enforcement $\vec{X}^{m+1,\ell+1}(\rho_{0,N})=\vec{r}(c_{l,r}^{m+1,\ell+1})$ during the iteration can destroy the inherent structure of the iterative algorithms. 
\end{rem}

 \begin{rem}
Our introduced methods relies on an exact arclength parameterization of the substrate curve, which is not always attainable for general substrates. In practice, we first discretize the substrate curve and compute the arclength at discrete points, and then adopt interpolation to solve its arclength parameterization. Nevertheless, in order to achieve a compatible interpolation error with the tolerance, the computational works can be rather heavy.  Thus it is important to find a smart approach to reduce the error when solving the arclength parameterization.
  \end{rem}

\subsection{Numerical results}

In this section, we present several numerical examples for SSD on various curved substrates. Unless otherwise stated, we always fix the material parameters $\eta=100$, $\sigma=-\frac{\sqrt{3}}{2}$. For the iterative solution of the structure-preserving scheme, we choose ${\rm tol}=10^{-9}$ in \eqref{eq:tol}. We restrict ourselves to the isotropic case and the anisotropic case as
\begin{itemize}
\item [(i)] the isotropic case with $\gamma(\vec{p})=|\vec{p}|$;
\item  [(ii)] the anisotropic case  with $l^4$-norm anisotropy which satisfies \eqref{eq:homoani}: 
\begin{equation}\label{eq:anipy}
\gamma(\vec{p}) =\gamma(p_1,p_2)= \sqrt[4]{(p_1)^4+(p_2)^4}.
\end{equation}
\end{itemize}
In the case of isotropy, we could simply choose $k(\vec{n}) = 2$ such that $\mathbf{Z}_k(\vec{n}) = I_2$, which implies \eqref{eq:newconvex} directly. For the anisotropy \eqref{eq:anipy}, we could set $k(\vec{n}) = 2\gamma(\vec{n})^{-3}$ so that \eqref{eq:newconvex}  is 
satisfied \cite{BJY21symmetrized}. Therefore in both cases, we have the condition \eqref{eq:zkassump} to be satisfied and thus the unconditional stability in \eqref{eq:denergylaw}. 

\begin{figure}[!htp]
  \centering
  \includegraphics[width=0.45\textwidth]{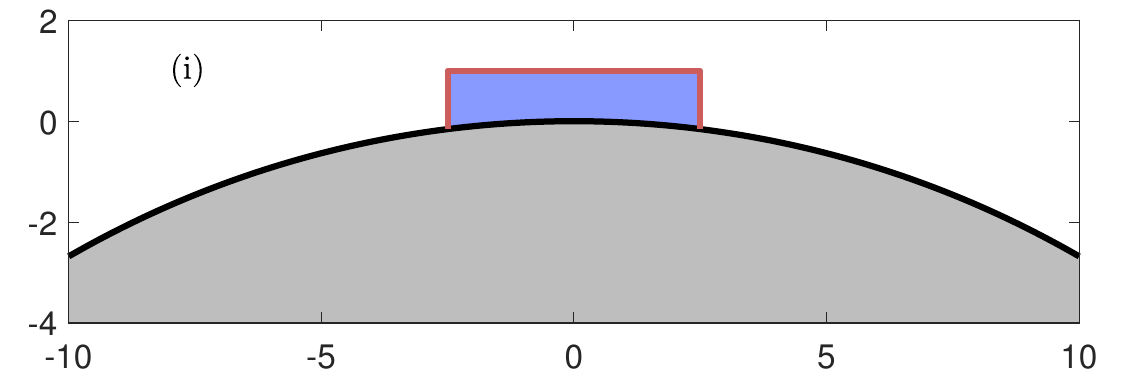}
  \includegraphics[width=0.45\textwidth]{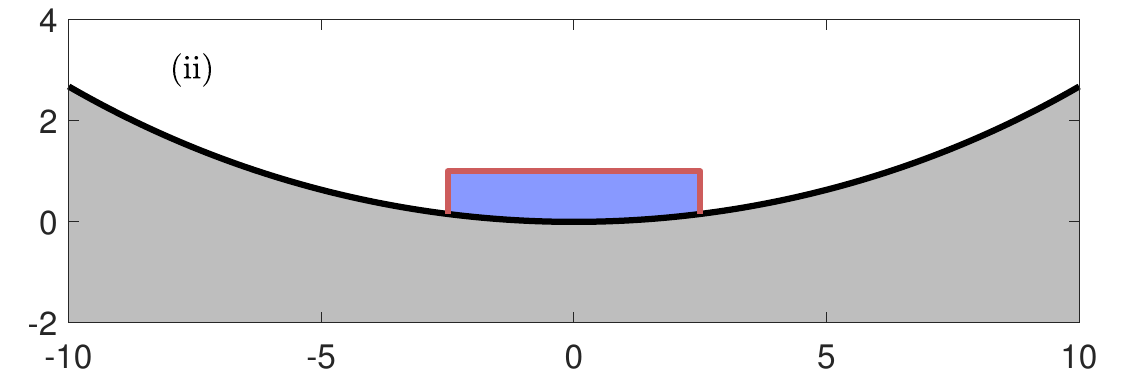}\\[0.7em]
  \includegraphics[width=0.45\textwidth]{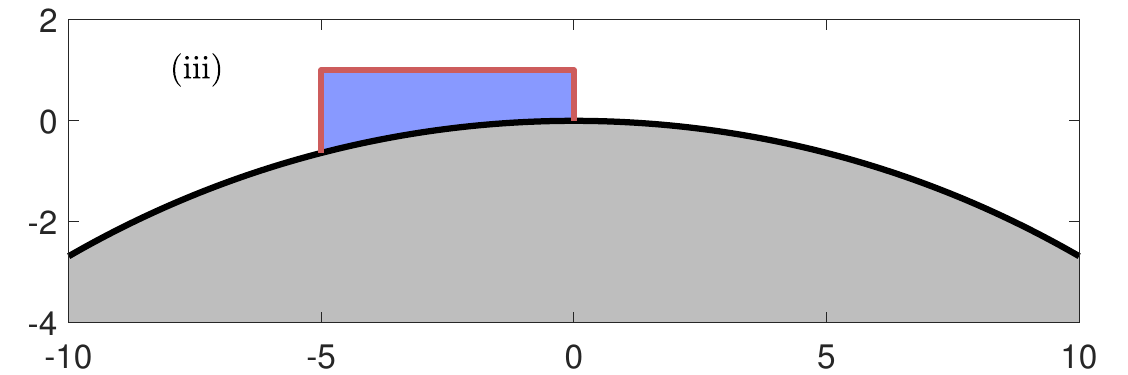}\includegraphics[width=0.45\textwidth]{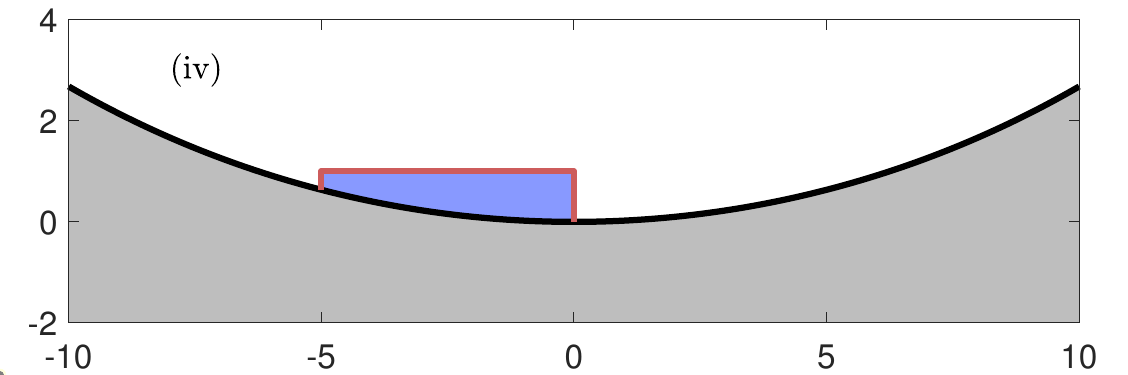}
  \caption{Four different initial setups for a initially deposited thin film on convex$\slash$concave substrates with different positions (i): a symmetric film on a concave substrate; (ii) a symmetric film on a convex substrate; (iii) a non-symmetric film on a concave substrate; (iv) a non-symmetric film on a convex substrate. }
  \label{fig:setup}
  \end{figure} 
  
\vspace{0.4em}
\noindent
{\bf Example 1}: We first test the convergence of our introduced method \eqref{eqn:vfd} in both isotropic case and anisotropic case with anisotropy given by \eqref{eq:anipy}. Initially, we choose a thin film of size $5\times 1$ and consider the following four different initial setting, as shown in Fig.~\ref{fig:setup}. The four cases are denoted by case (i) on a concave substrate with symmetry; case (ii) on a convex substrate with symmetry; case (iii) on a concave substrate with asymmetry; case (iv) on a convex substrate with asymmetry. The substrate curve is given by  a circle of radius $R=20$.

Given the numerical solutions $\{\vec{X}^m\}_{0\leq m\leq M}$ using mesh size $h$ and $\ttau$, we define the numerical approximation solution 
\[\vec{X}_{h, \ttau}(t, \rho_j):= \frac{t-t_{m}}{\ttau}\vec{X}^{m}(\rho_j)+ \frac{t_{m+1}-t}{\ttau}\vec{X}^{m+1}(\rho_j), \quad \forall j = 1, 2, \ldots, N-1,\]
with boundary given by:
\begin{equation*}
\vec{X}_{h, \ttau}(t, \rho_0):=\vec{r}\left(\frac{t-t_{m}}{\ttau}c_l^{m}+ \frac{t_{m+1}-t}{\ttau}c_l^{m+1}\right), \quad
\vec{X}_{h, \ttau}(t, \rho_N):=\vec{r}\left(\frac{t-t_{m}}{\ttau}c_r^{m}+ \frac{t_{m+1}-t}{\ttau}c_r^{m+1}\right).
\end{equation*}
 We then measure the numerical error $e^{h,\ttau}(t)$ by comparing $\vec{X}_{h,\ttau}(t)$ with the reference solution $\vec{X}_r(t)$, which is obtained with a refined mesh size $h=2^{-8}$ and time step $\ttau = 2^{-16}$. The measurement is calculated using the manifold distance discussed in \cite{Zhao20}, i.e., the area of the symmetric difference region:
 \[e^{h,\ttau}(t):= {\rm MD}\left(\vec{X}_{h,\ttau}(t), \vec{X}_r(t)\right)=|(\Omega_{h,\ttau}(t)\backslash \Omega_r(t))\cup(\Omega_r(t)\backslash\Omega_{h,\ttau}(t))|, \]
where $\Omega_{h,\ttau}(t)$ and $\Omega_r(t)$ represents the area of the region occupied by the thin film with interface $\vec{X}_{h,\ttau}(t)$ and $\vec{X}_r(t)$, respectively, and $|\Omega|$ denote the area of the region $\Omega$.

\begin{figure}[t]
\centering
 \includegraphics[width=0.7\textwidth]{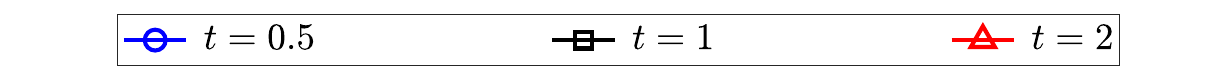}
\includegraphics[width=0.45\textwidth]{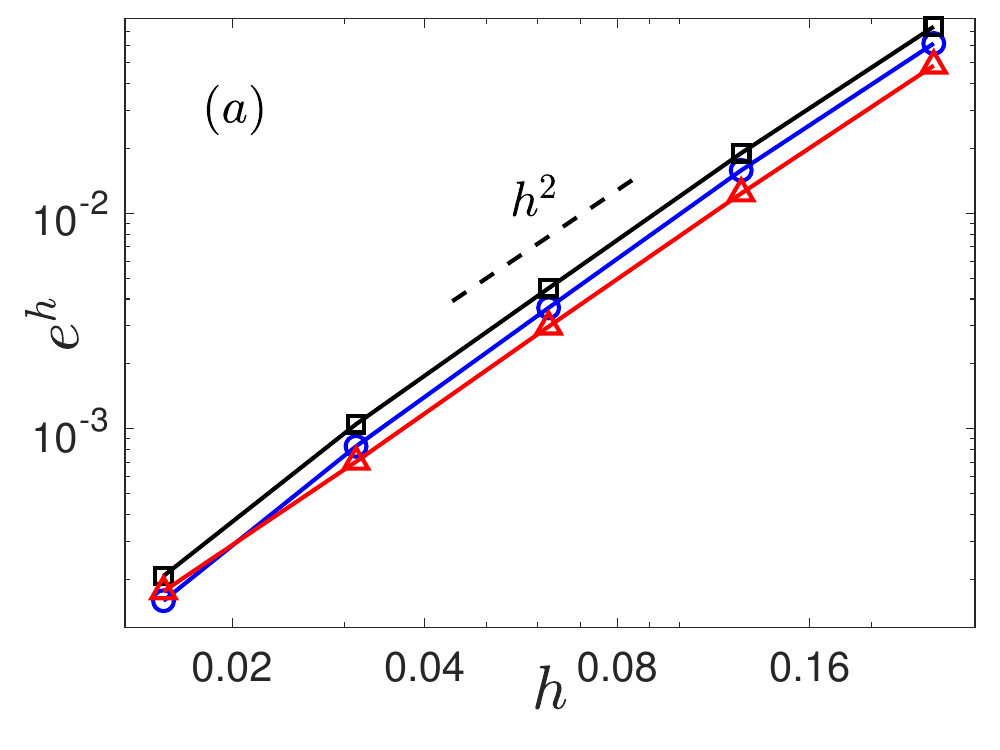}\hspace{0.7em}\includegraphics[width=0.45\textwidth]{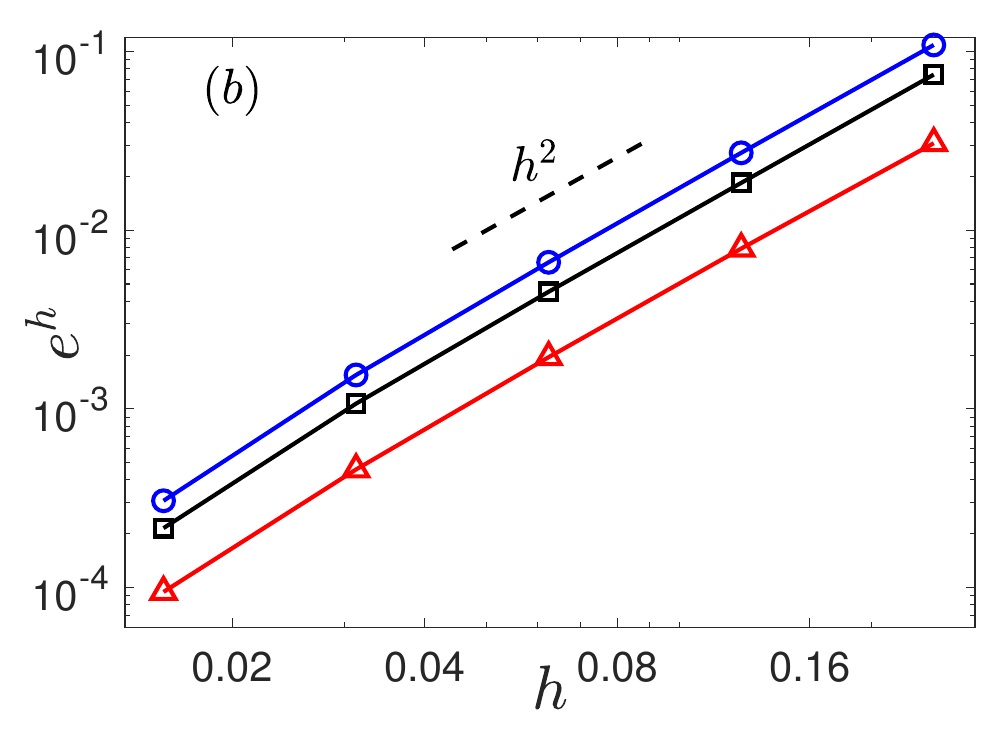}
\includegraphics[width=0.45\textwidth]{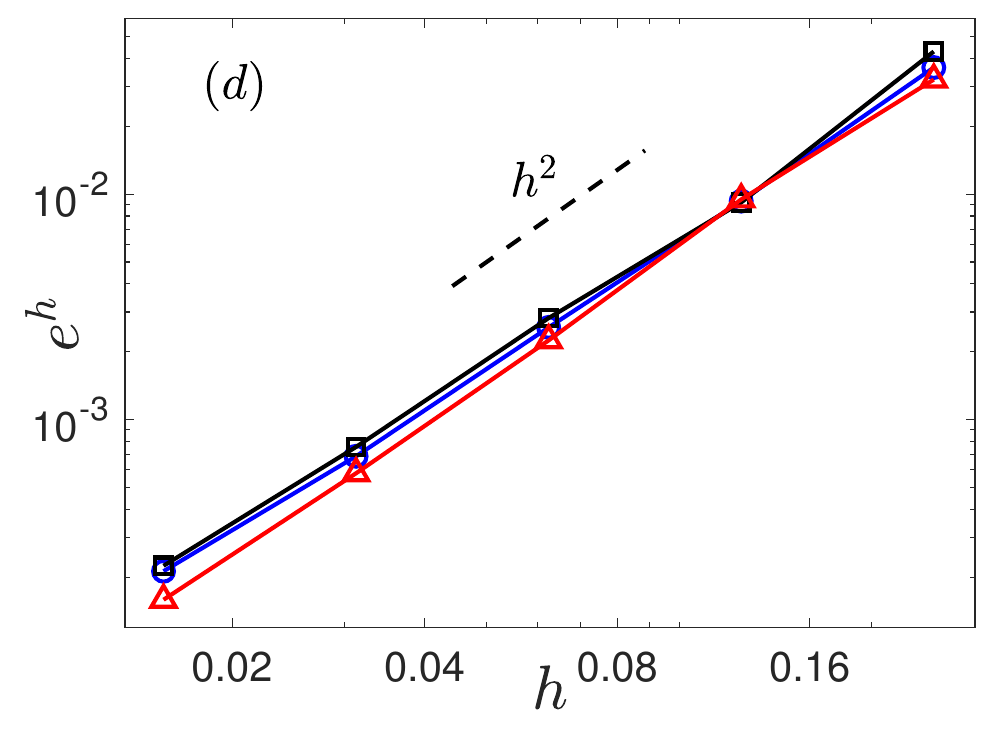}\hspace{0.7em}\includegraphics[width=0.45\textwidth]{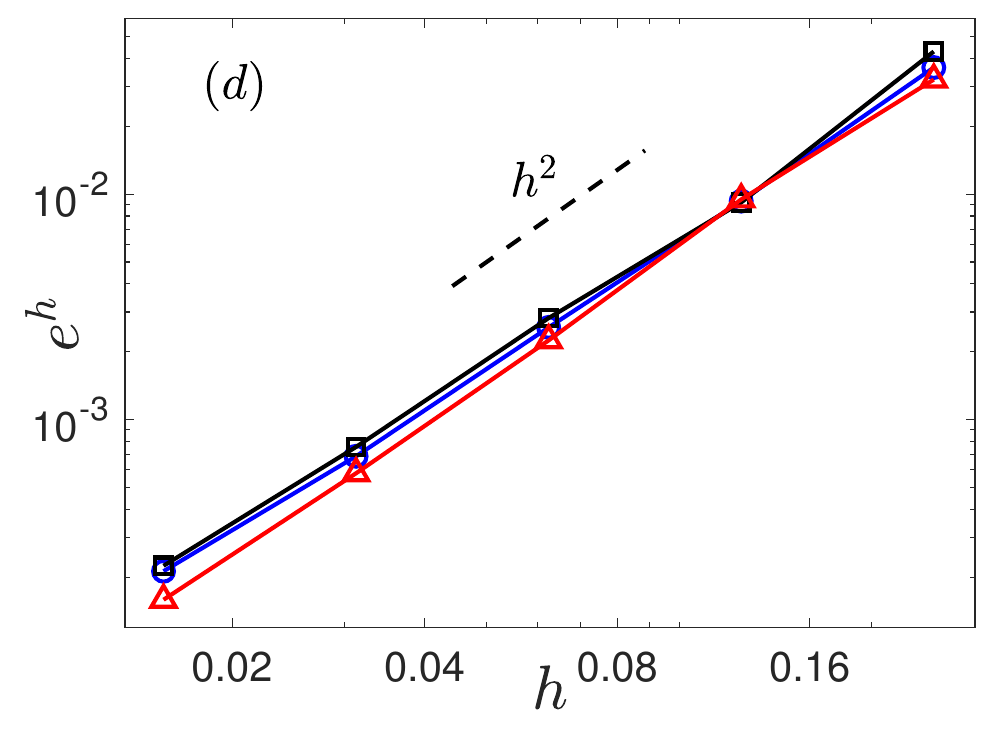}
\caption{Convergence test with different setups : (a) case (i) with isotropy ; (b) case (ii) with isotropy; (c) case (iii) with anisotropy; (d) case (iv) with anisotropy. }
\label{fig:convergence}
\end{figure}

We fix $\ttau = h^2$ and the convergence results are shown in Fig.~\ref{fig:convergence} for the four difference cases. Here we observe a second-order convergence rate with respect to the mesh size $h$ for different setups, which demonstrates the robustness and efficiency of the introduced method.

\begin{figure}[!htp][b]
  \centering
  \includegraphics[width=0.7\textwidth]{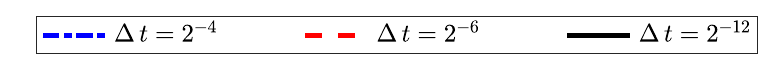}
  \includegraphics[width=0.45\textwidth]{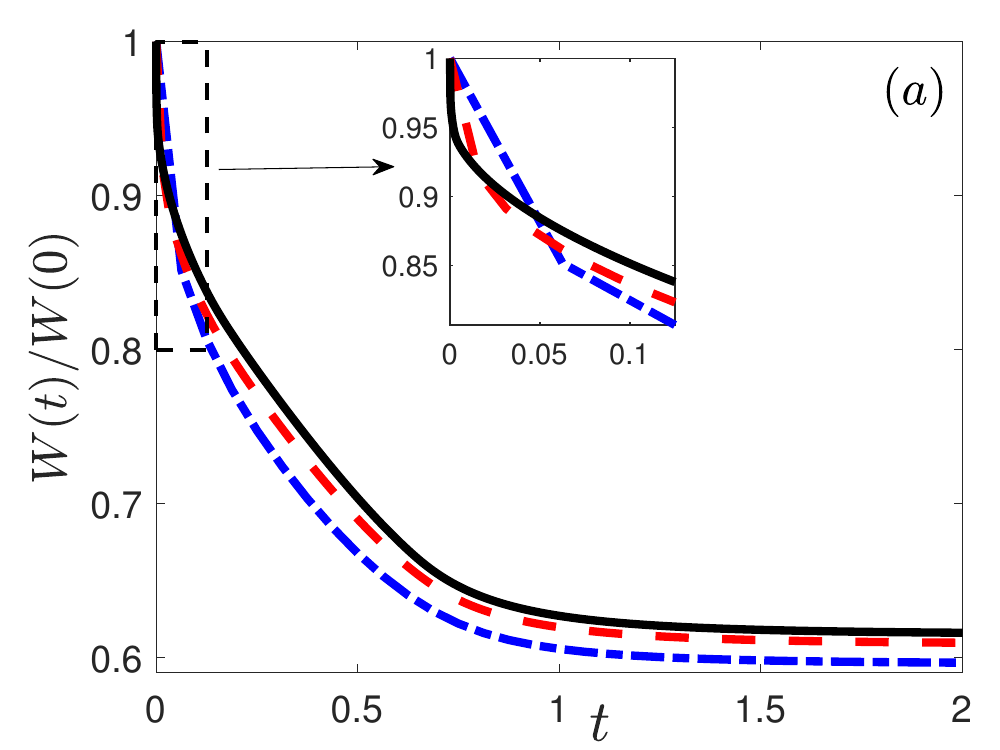}\hspace{0.7em}
  \includegraphics[width=0.45\textwidth]{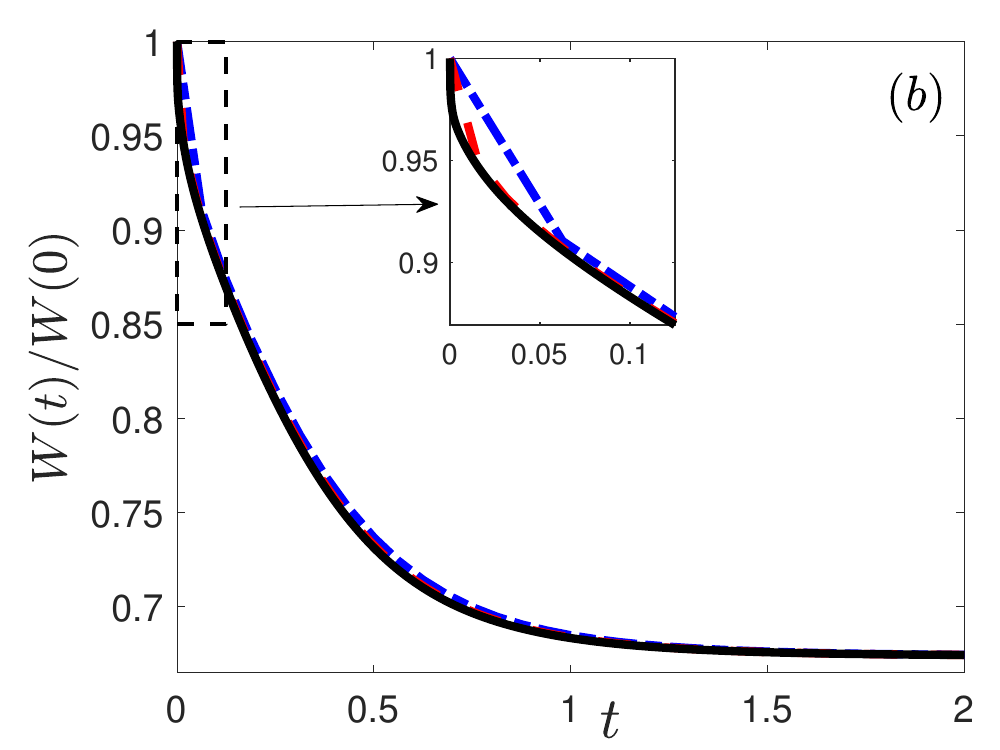}
  \includegraphics[width=0.45\textwidth]{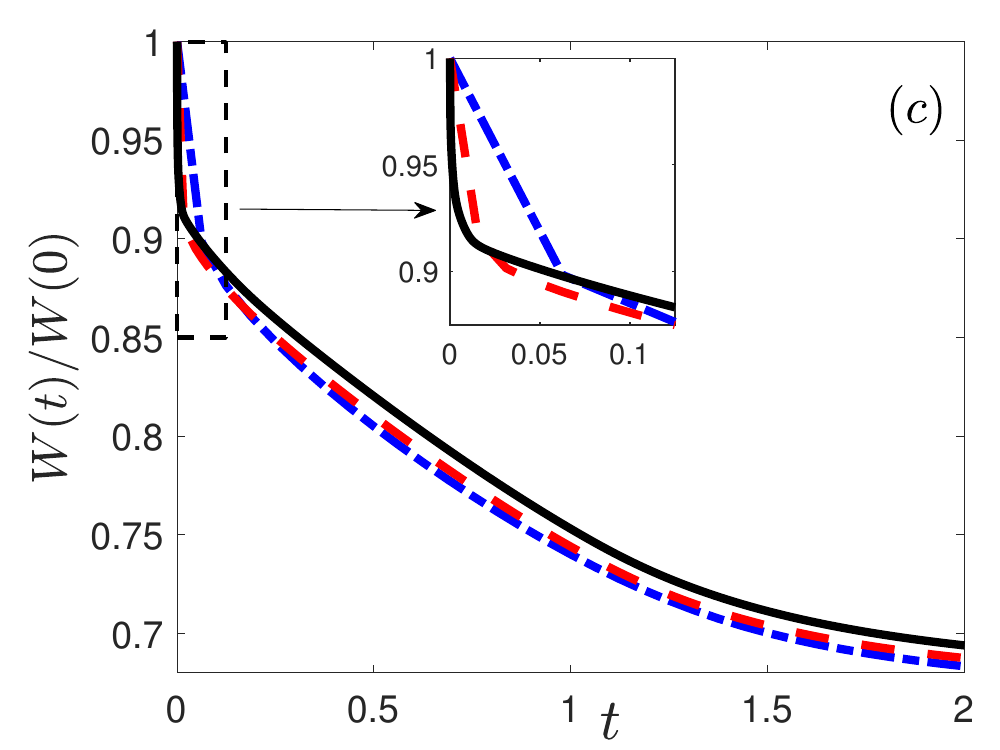}\hspace{0.7em}
  \includegraphics[width=0.45\textwidth]{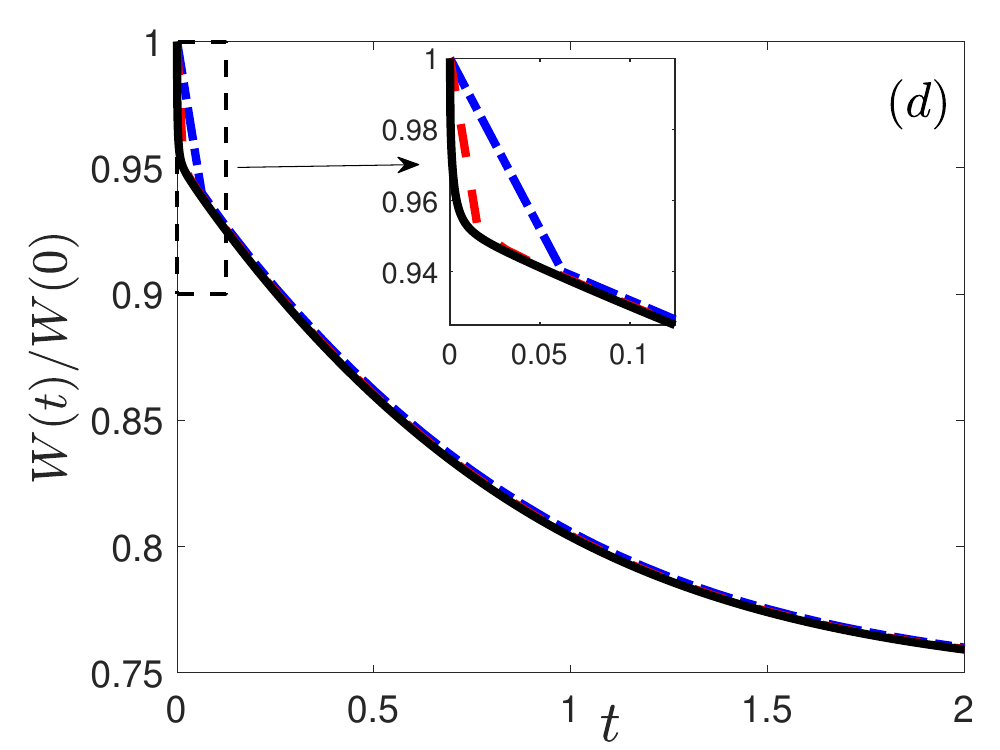}
  \caption{The time history plots of the normalized energy under different time steps  (a) case (i) with isotropy; (b) case (ii) with isotropy; (c) case (iii) with anisotropy; (d) case (iv) with anisotropy. }
  \label{fig:outer}
  \end{figure}
  
\begin{figure}[!htp]
    \centering
    \includegraphics[width=0.8\textwidth]{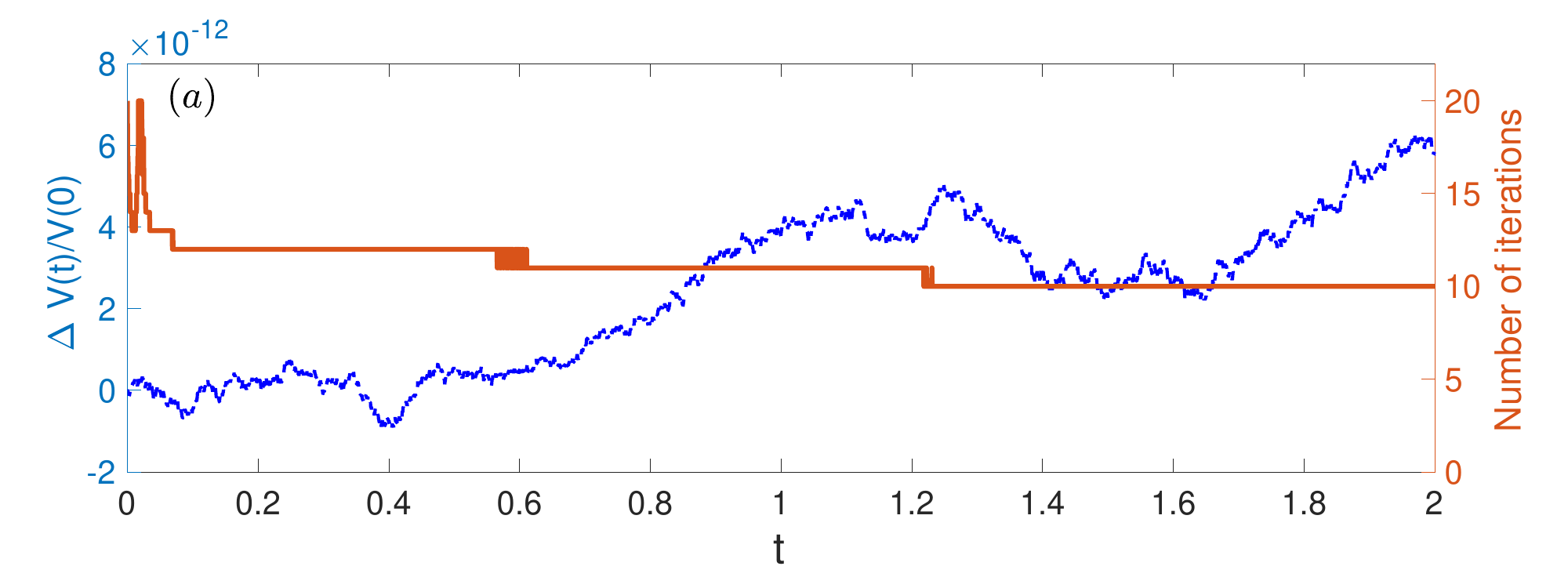}
    \includegraphics[width=0.8\textwidth]{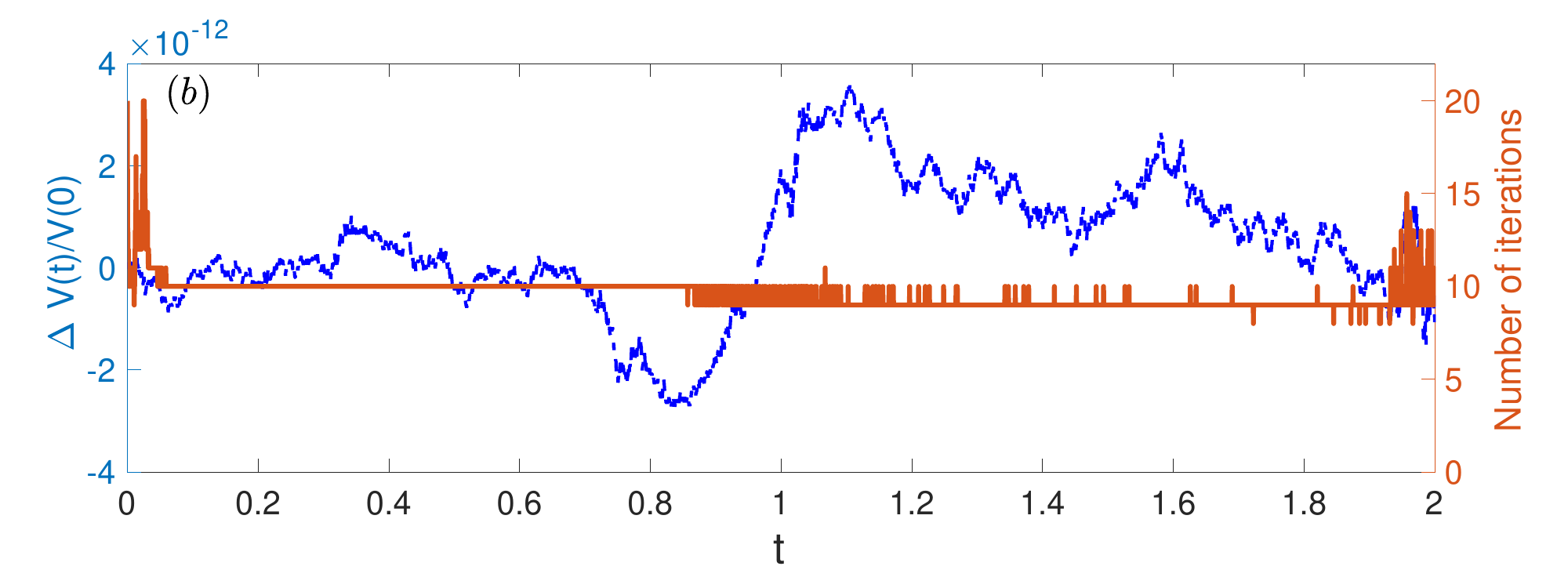}
    \caption{Mass preservation and the number of iterations with $\ttau = 2^{-14}$ and $h=2^{-7}$: (a) case (i) with isotropy; (b) case (ii) with anisotropy.}
    \label{fig:area}
    \end{figure}

\vspace{0.4em}
\noindent
{\bf  Example 2}: To further illustrate the structure-preserving properties of our introduced method, we use the same setup as in  Example 1. The numerical results for the energy decay and area preservation are presented in Fig.~\ref{fig:outer} and Fig.~\ref{fig:area}, respectively. We observe that a monotonic decay of the discrete energy even for  pretty large time steps, and the mass is well preserved up to $10^{-12}$. Moreover, for the introduced iterative algorithm in \S\ref{sec:sol}, we observe mostly only 10-12 number of iterations are required.  These numerical results confirm our theoretic results in Theorem \ref{thm:sp} as well as the computational efficiency of the solution method for the nonlinear method. 

\begin{figure}[!htp]
\centering
\includegraphics[width=0.45\textwidth]{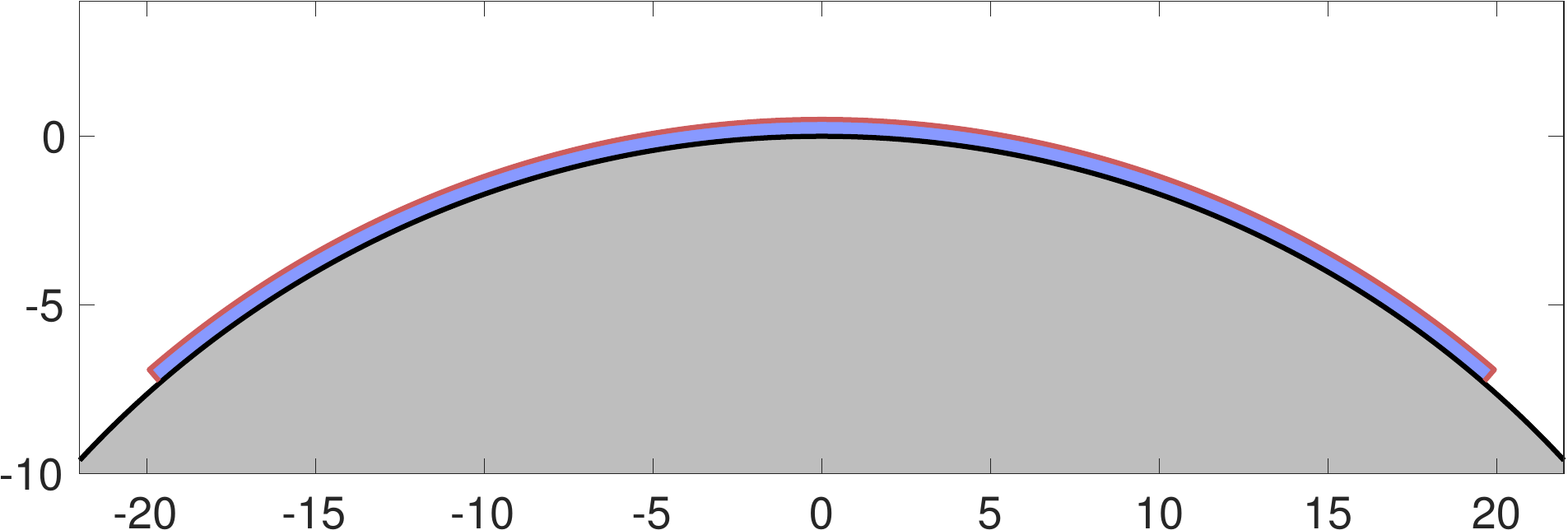}\hspace{0.7em}\includegraphics[width=0.45\textwidth]{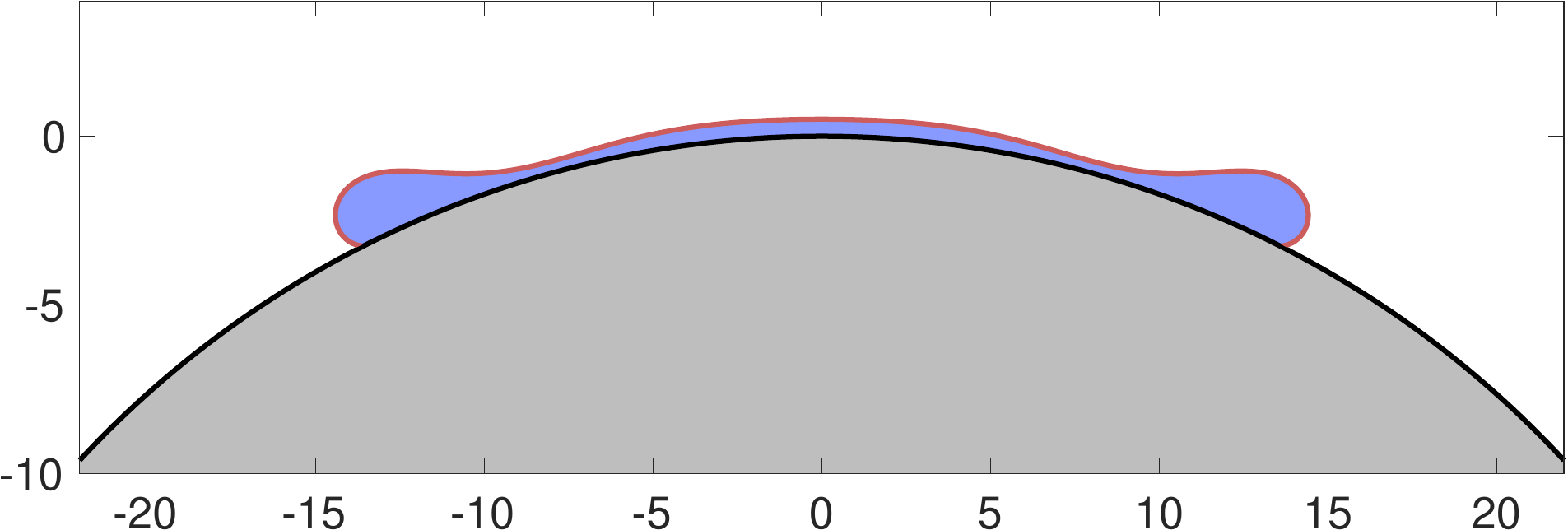}\\[0.8em]
\includegraphics[width=0.45\textwidth]{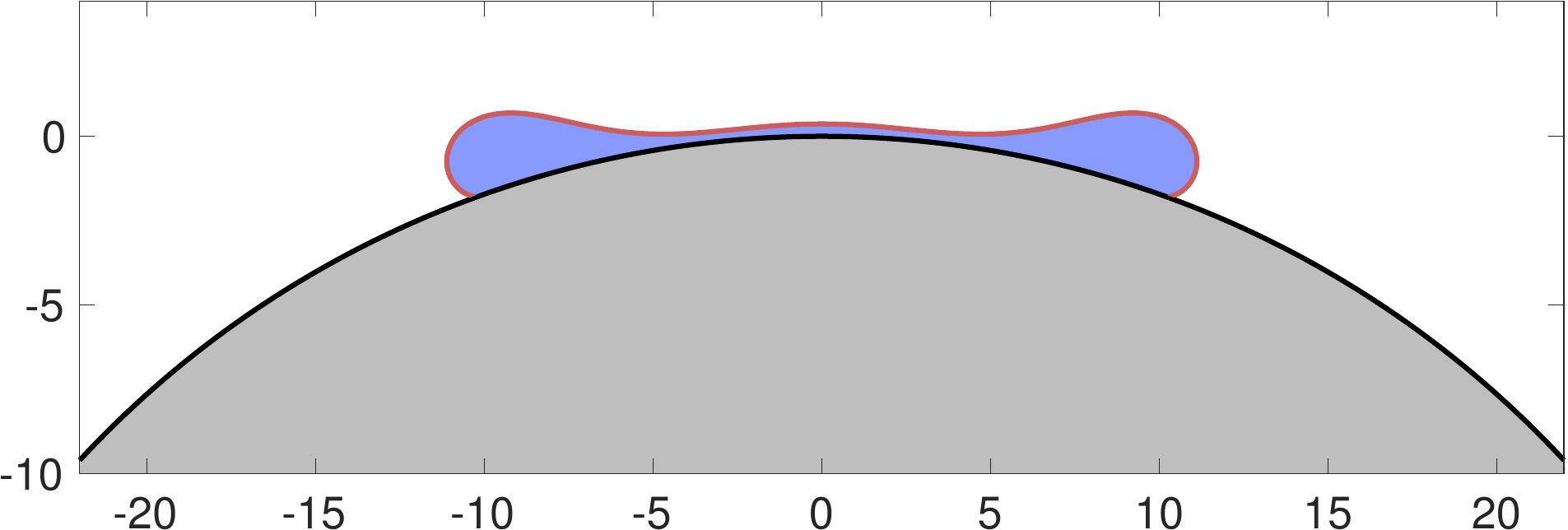}\hspace{0.7em}\includegraphics[width=0.45\textwidth]{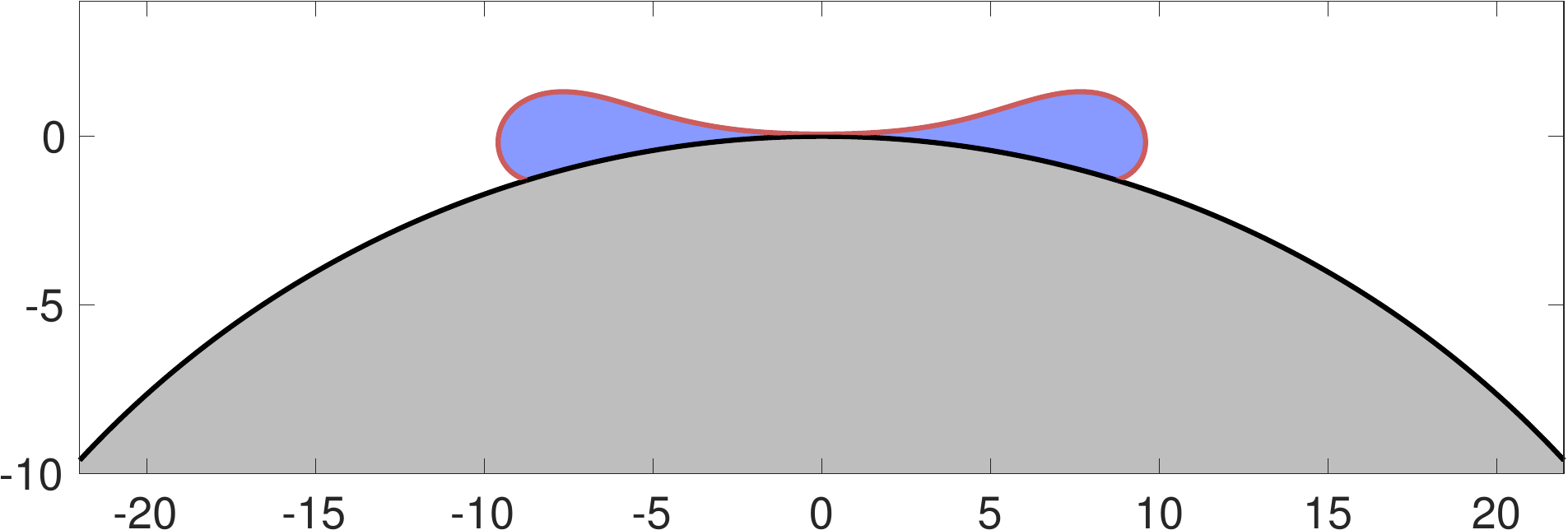}\\[0.7em]
\includegraphics[width=0.45\textwidth]{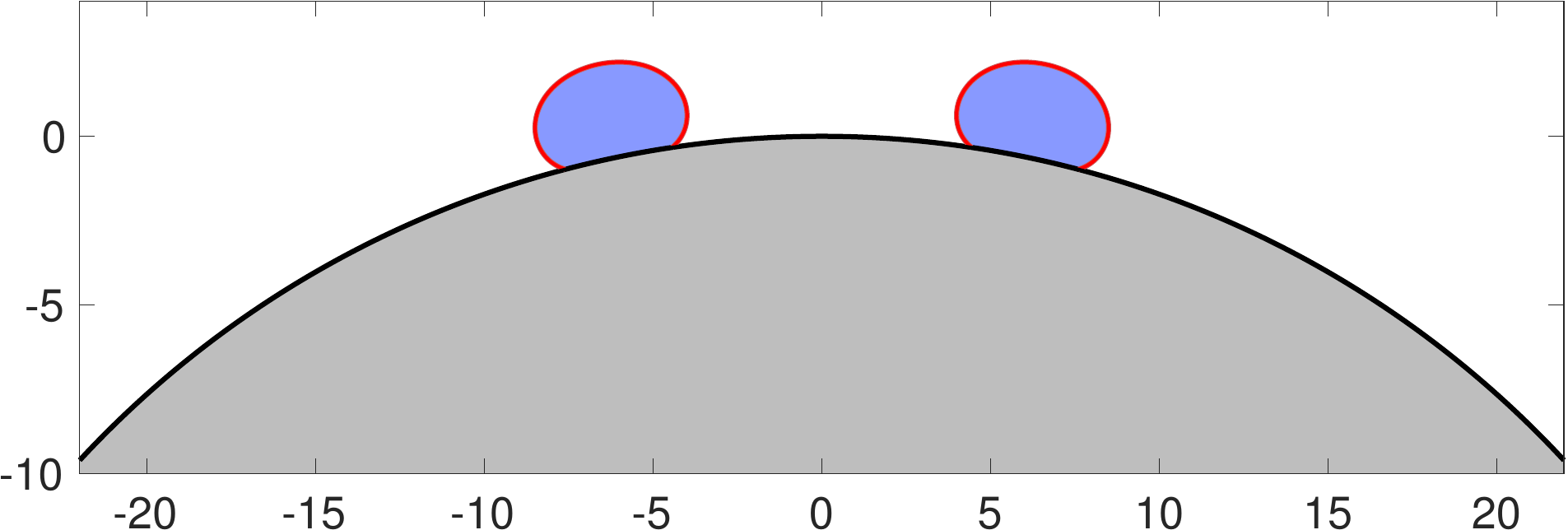}\hspace{0.8em}\includegraphics[width=0.45\textwidth]{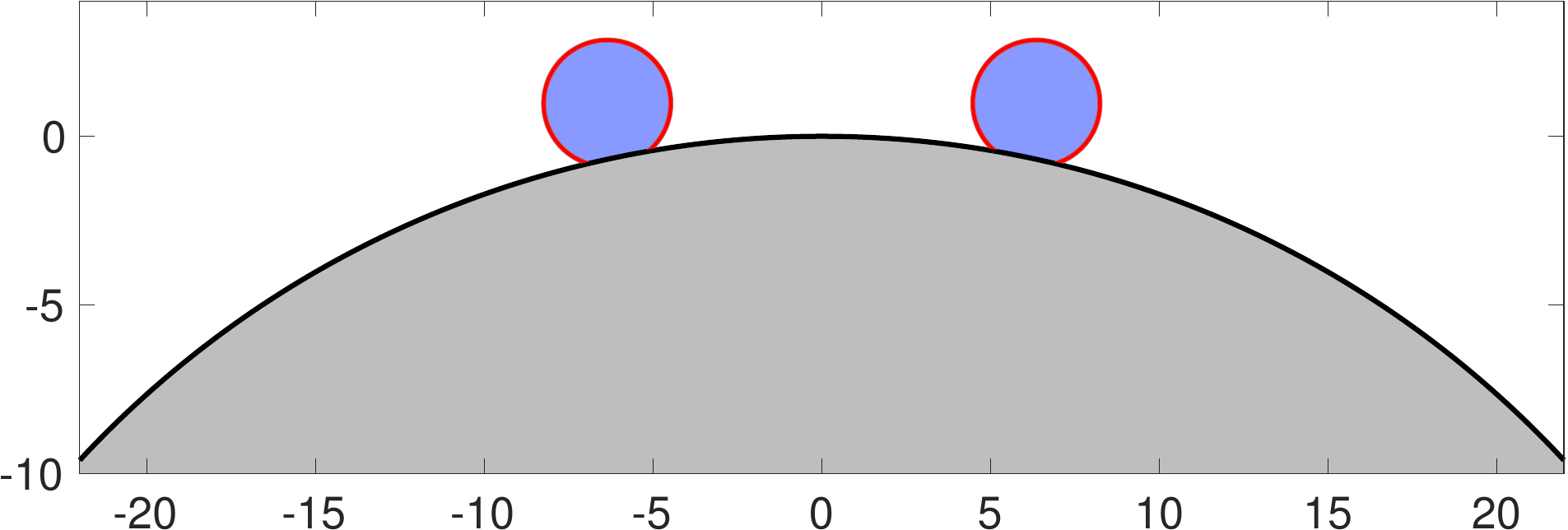}
\caption{Snapshots in the evolution of a long thin film on a convex substrate of circular shape at times $t=0, 15.625, 40.5, 56.25, 68.75, 93.75$ (from left to right and then from top to bottom). }
\label{fig:evolution_outer}
\end{figure}

\begin{figure}[!htp]
  \centering
  \includegraphics[width=0.45\textwidth]{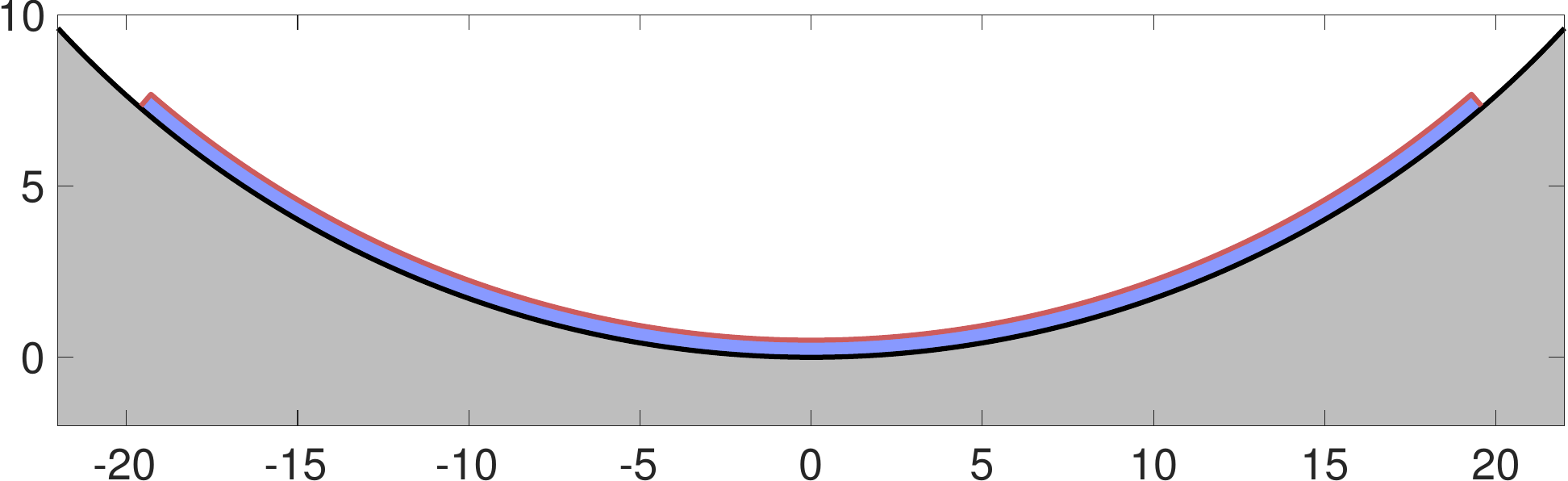} \hspace{0.7em}\includegraphics[width=0.45\textwidth]{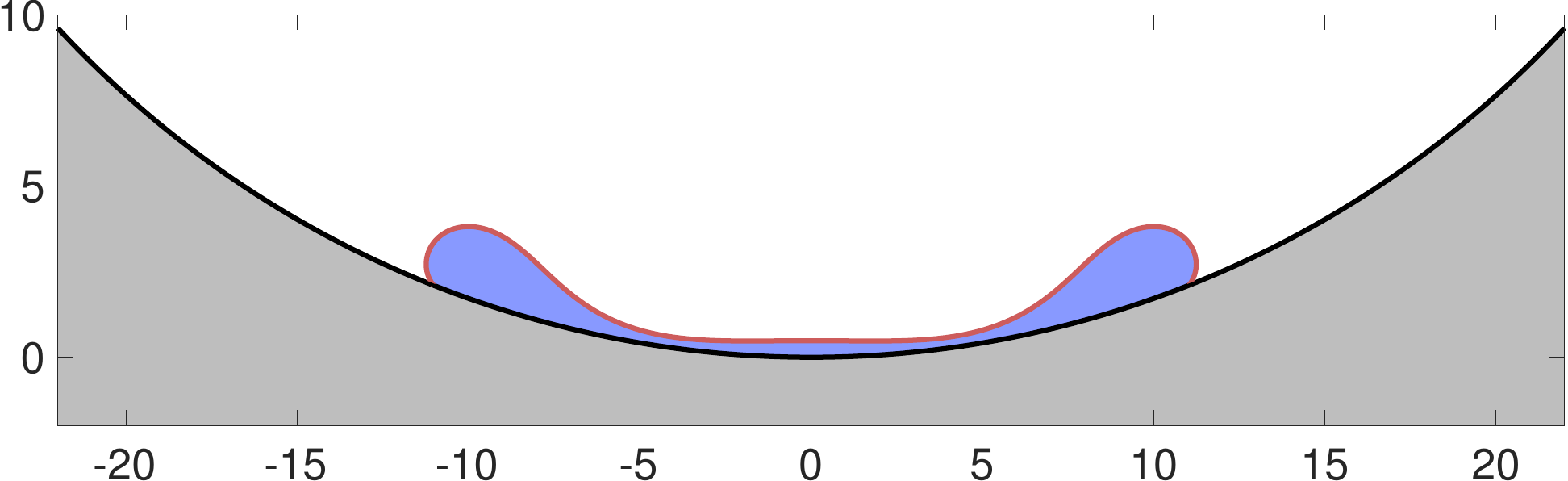}\\[0.7em]
  \includegraphics[width=0.45\textwidth]{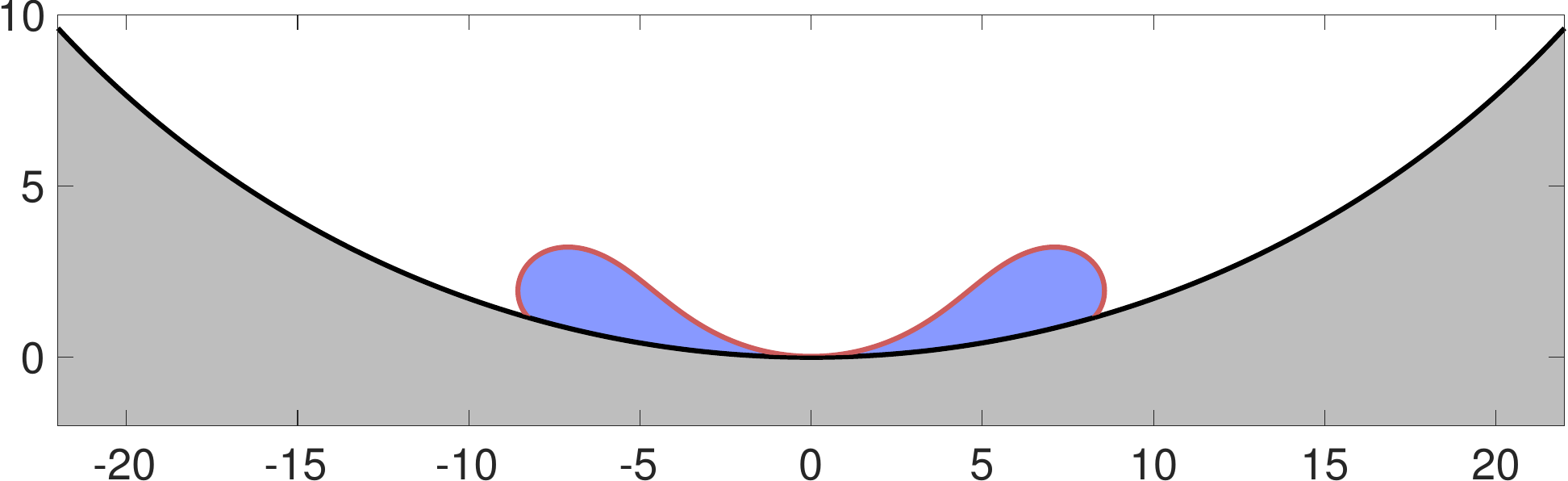} \hspace{0.7em}\includegraphics[width=0.45\textwidth]{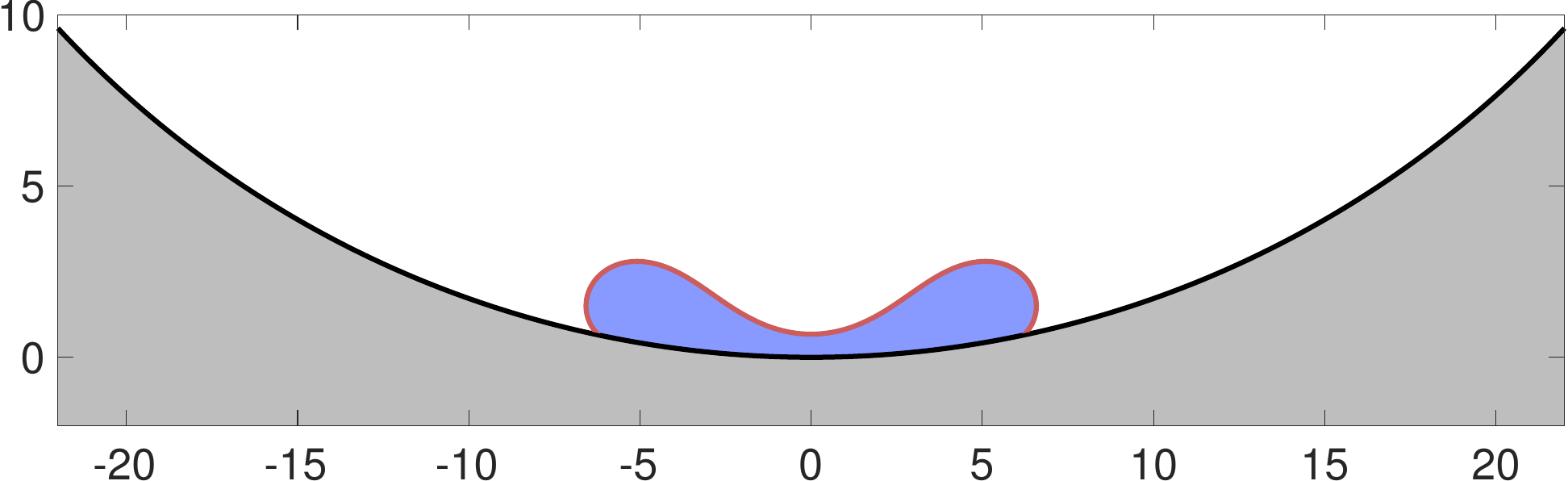}\\[0.7em]
  \includegraphics[width=0.45\textwidth]{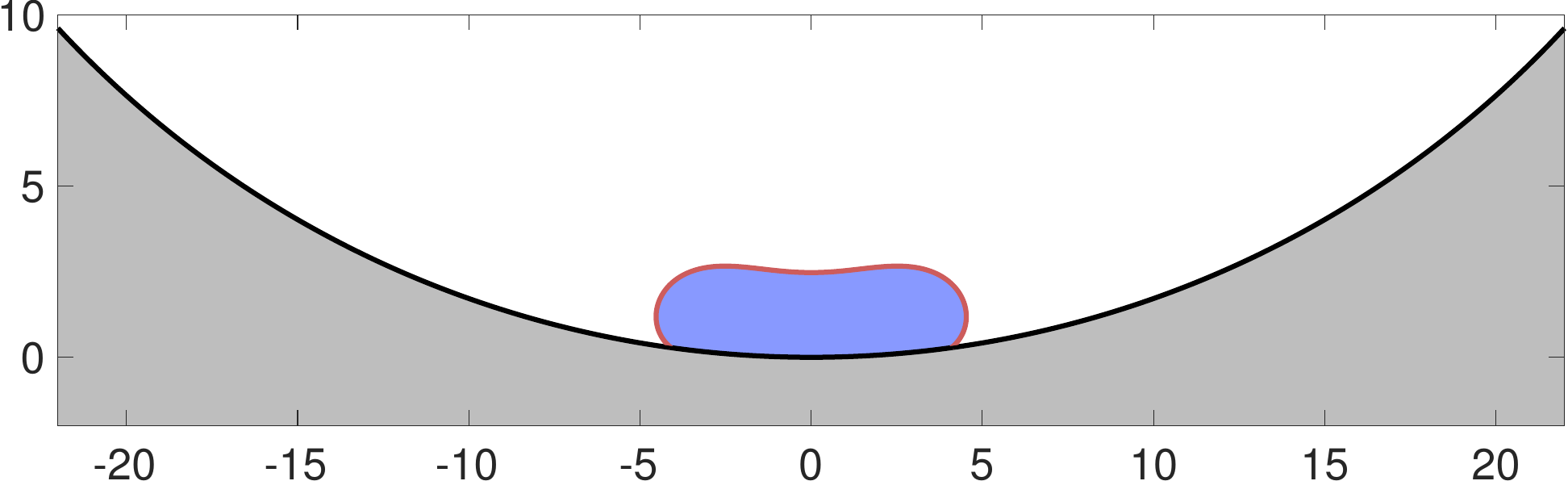} \hspace{0.7em}\includegraphics[width=0.45\textwidth]{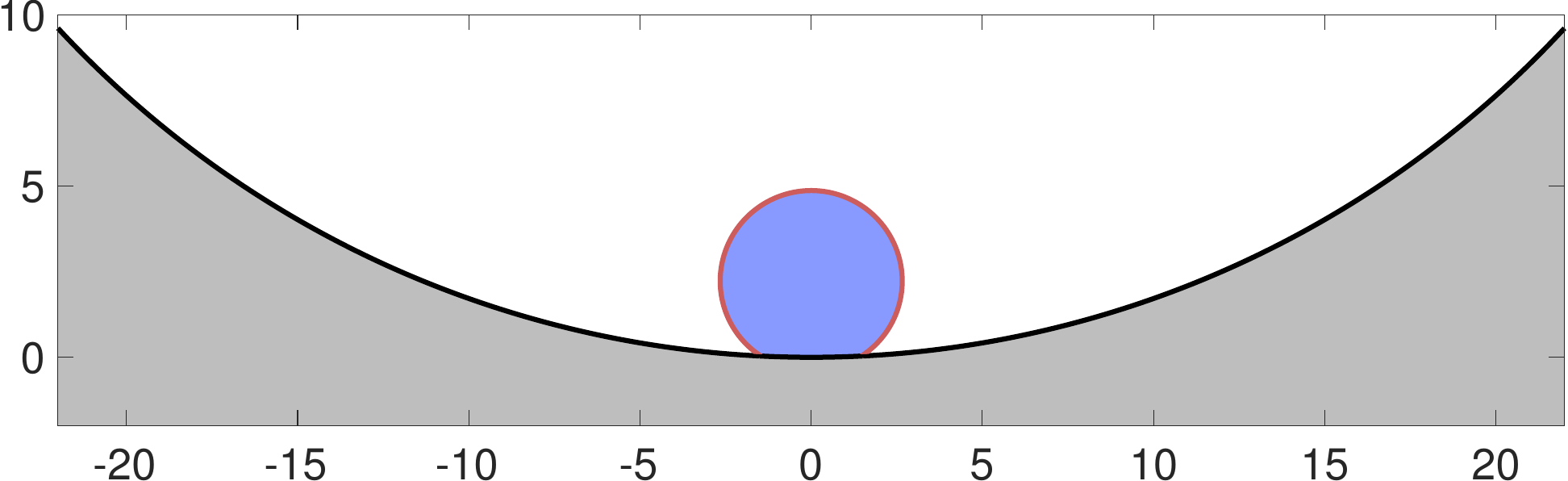}
  \caption{Snapshots in the evolution of a long thin film on a concave substrate of circular shape at times $t=0, 28.5, 52, 75, 93.75, 187.5$ (from left to right and then from top to bottom). }
  \label{fig:evolution_inner}
  \end{figure}

\vspace{0.4em}
\noindent
{\bf Example 3}: In this example, we consider the evolution of thin film on both convex and concave substrates, and the substrate curve is given by a circle of radius $R=30$. Initially the size of the thin film is aligned with the substrate curve with thickness $0.5$ and length $42$. Several snapshots of the thin film are depicted in Fig.~\ref{fig:evolution_outer} and Fig.~\ref{fig:evolution_inner} for the two different cases. Here we observe that the thin film breaks into two small islands on the convex substrate. On the contrary, the thin film on the concave substrate evolves into a single island as the equilibrium. This implies that islands on convex substrate tend to break up more easily than that on concave substrate. 

\begin{figure}[!htp]
  \centering
  \includegraphics[width=0.42\textwidth]{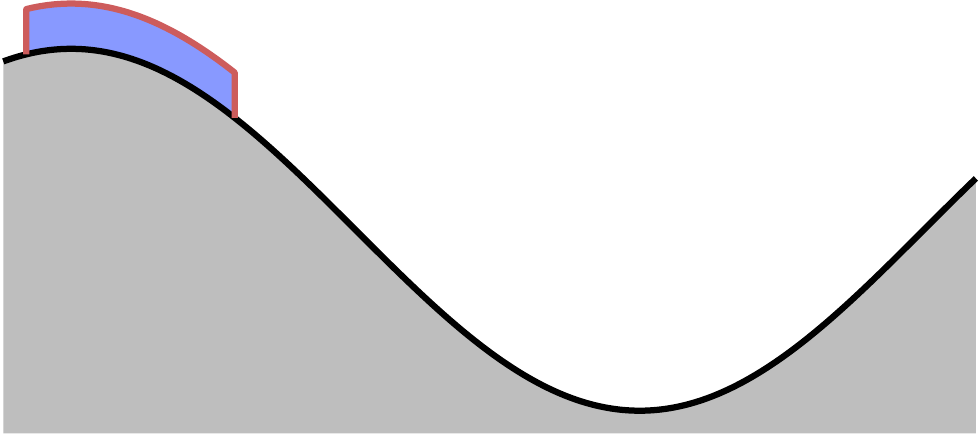}\hspace{0.9em}\includegraphics[width=0.42\textwidth]{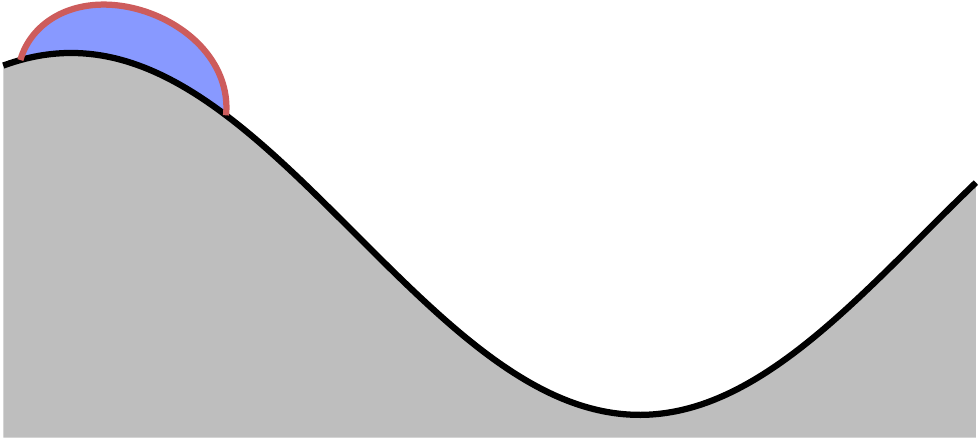}\\[0.8em]
  \includegraphics[width=0.42\textwidth]{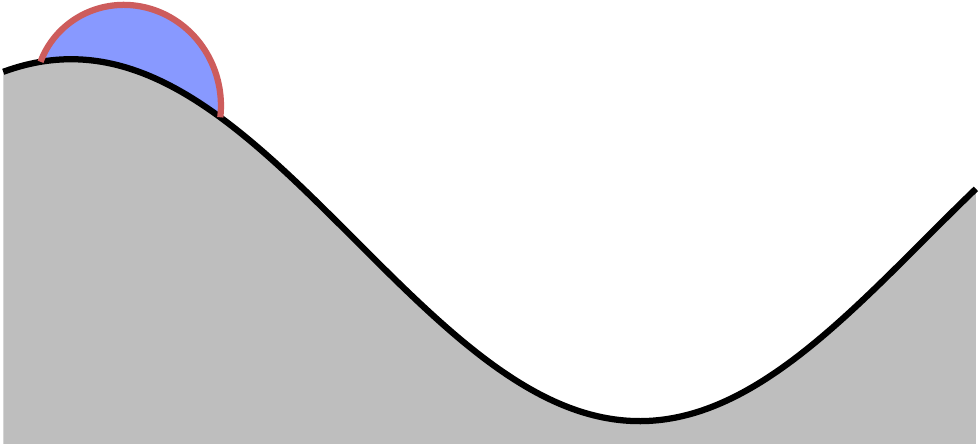}\hspace{0.9em}\includegraphics[width=0.42\textwidth]{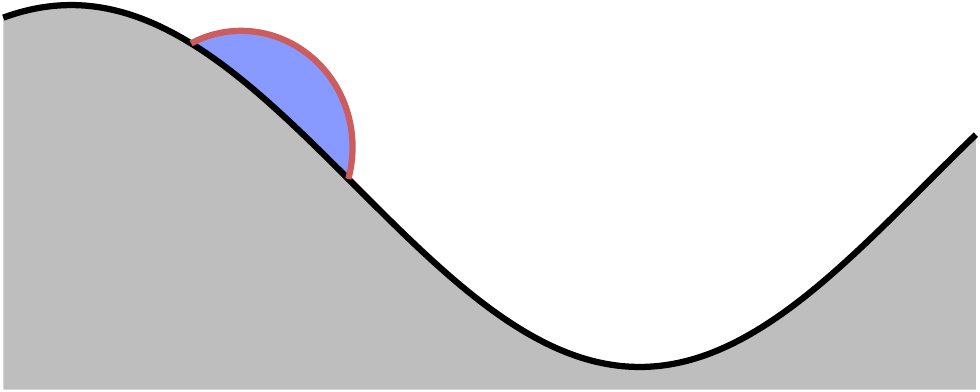}\\[0.8em]
  \includegraphics[width=0.42\textwidth]{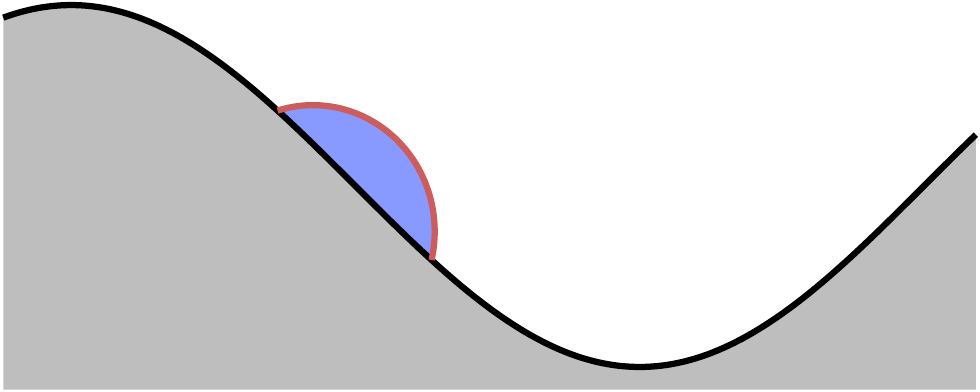}\hspace{0.9em}\includegraphics[width=0.42\textwidth]{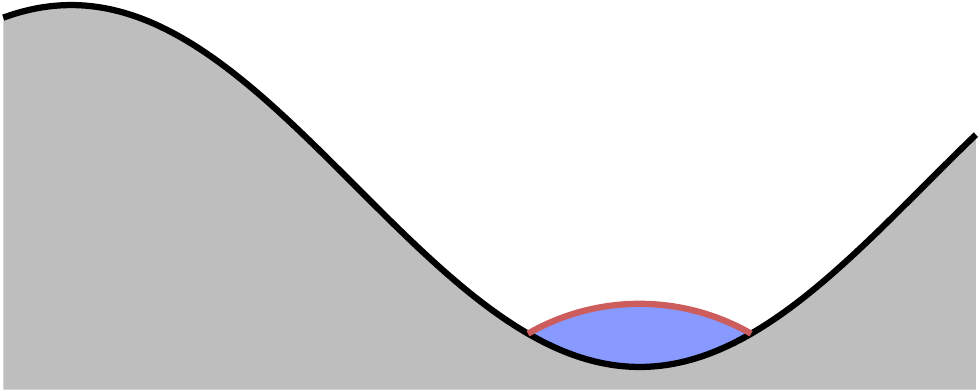}
  \caption{Snapshots in the evolution of a short thin film on a  sinusoidal substrate curve with isotropy at times $t=0, 1.25, 12.5, 250, 500, 1225$ (from left to right and then from top to bottom).}
  \label{fig:evolution_isoshort}
  \end{figure}

\begin{figure}[!htp]
    \centering
    \includegraphics[width=0.4\textwidth]{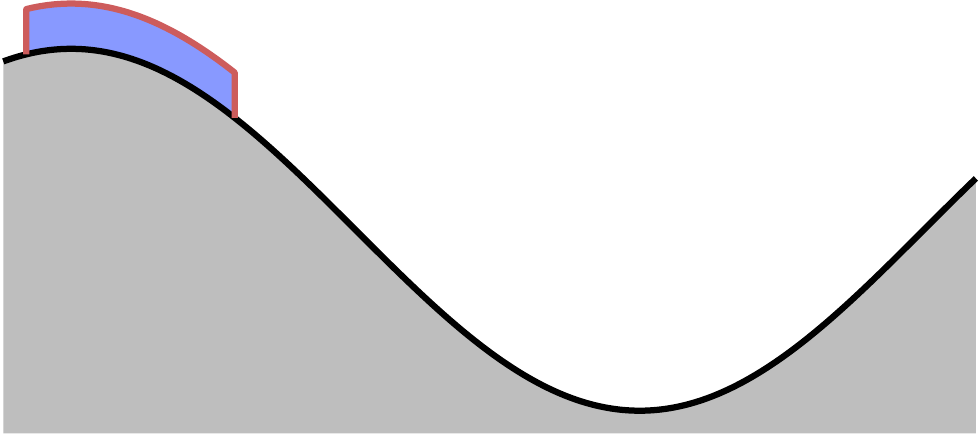}\hspace{0.9em}\includegraphics[width=0.4\textwidth]{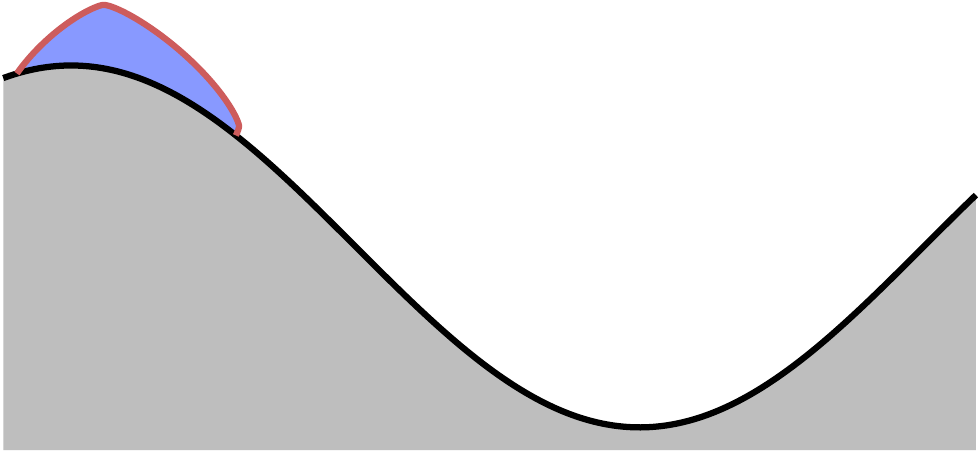}\\[0.8em]
    \includegraphics[width=0.4\textwidth]{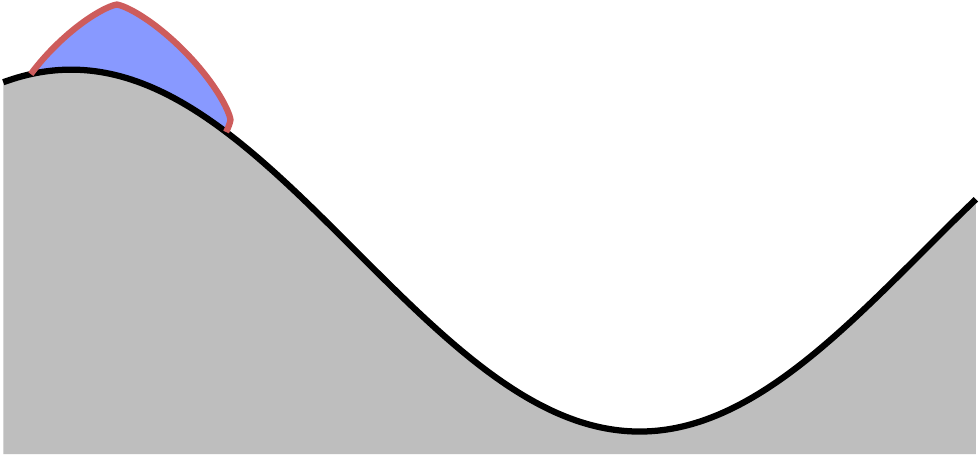}\hspace{0.9em}\includegraphics[width=0.4\textwidth]{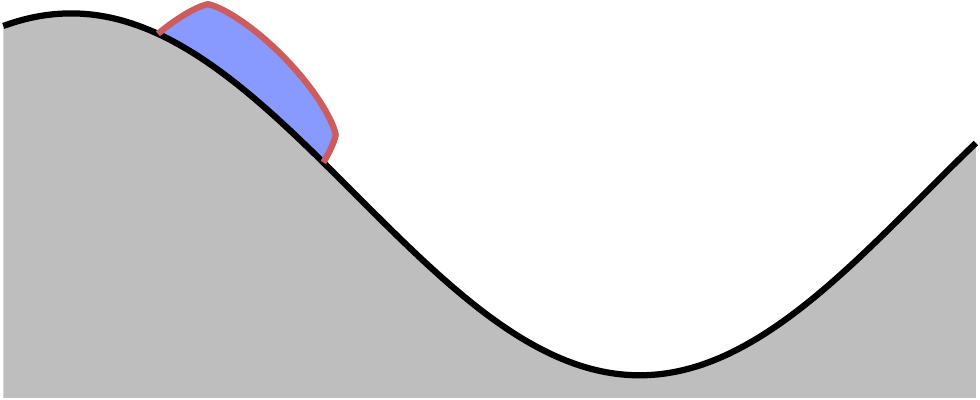}\\[0.8em]
    \includegraphics[width=0.4\textwidth]{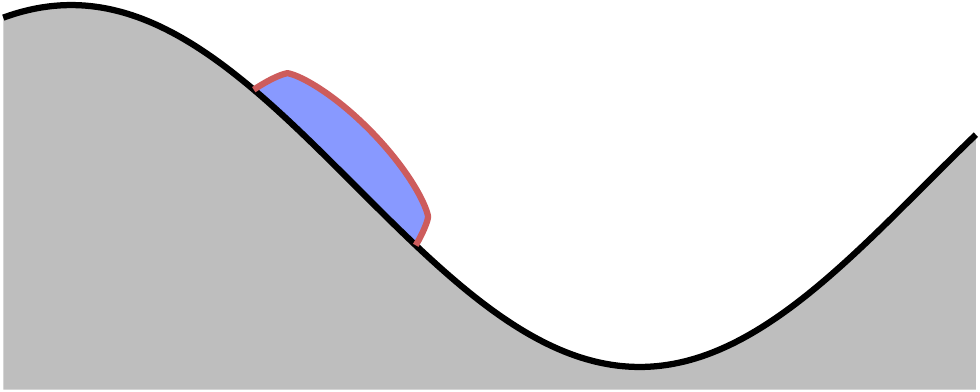}\hspace{0.9em}\includegraphics[width=0.4\textwidth]{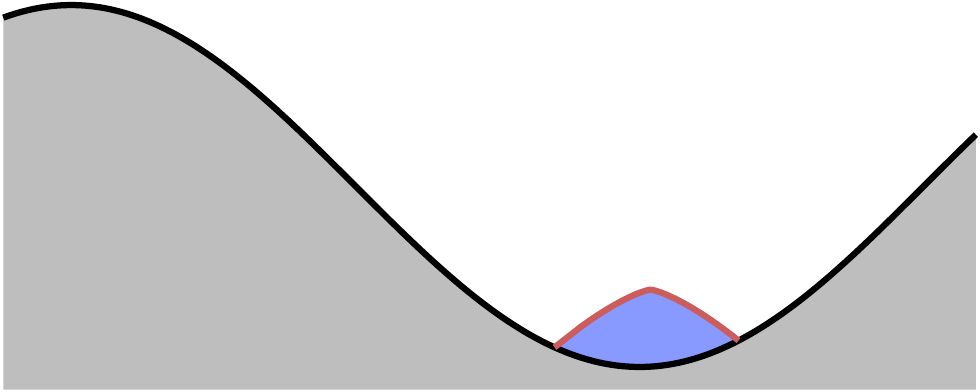}
    \caption{Snapshots in the evolution of a short thin film on a sinusoidal substrate curve with anisotropy at times $t=0, 1.25, 12.5, 250, 500, 1125$ (from left to right and then from top to bottom).}
    \label{fig:evolution_anisoshort}
    \end{figure}

\begin{figure}[!htp]
  \centering
  \begin{tabular}{cc}
    \includegraphics[width=0.4\textwidth]{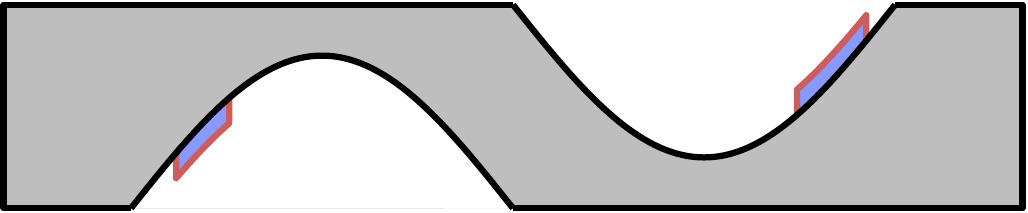} &
    \includegraphics[width=0.4\textwidth]{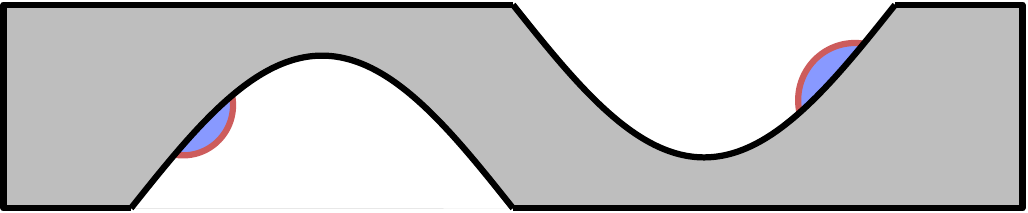} \\[0.7em]
    \includegraphics[width=0.4\textwidth]{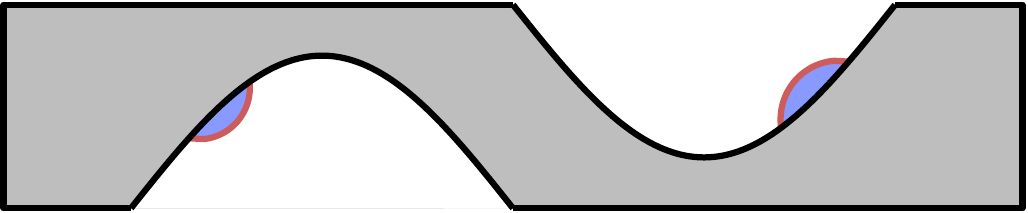} &
    \includegraphics[width=0.4\textwidth]{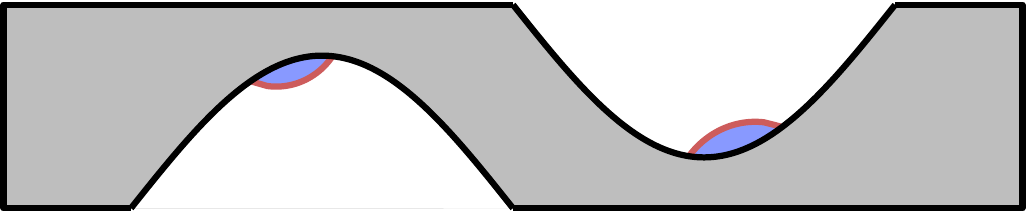} \\[0.7em]
    \includegraphics[width=0.4\textwidth]{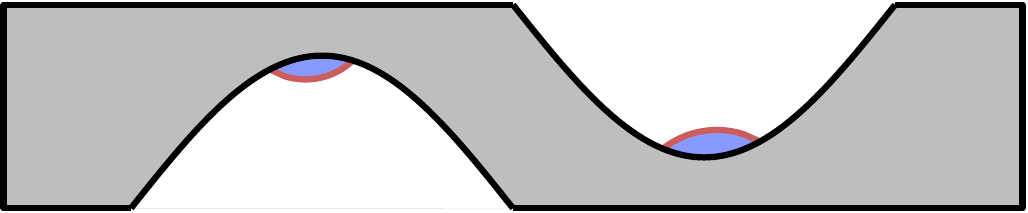} &
    \includegraphics[width=0.4\textwidth]{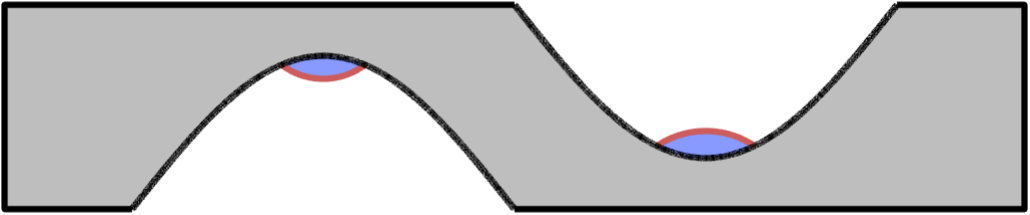} 
  \end{tabular}
  \caption{Snapshots in the evolution of two small thin films on a periodically closed sinusoidal substrate curve with isotropy at times $t=0, 0.156, 0.391, 1.562,3.906, 15.625$ (from left to right and then from top to bottom).}
  \label{fig:evolution_convexconcave}
  \end{figure}

\vspace{0.4em}
\noindent
{\bf Example 4}: We next examine the dewetting of thin films on general curve substrates, and the substrate curve is modeled by a sinusoidal shape given by $y=4\cos(x/4)$. Fig. \ref{fig:evolution_isoshort} and Fig. \ref{fig:evolution_anisoshort} show the snapshots of the thin films in isotropic and anisotropic cases, respectively.  As can be seen, the thin film migrates towards the lower curvature sites of the substrate, which shows high consistency to the results in \cite{ahn80,Klinger12,Jiang18curved, Zhao24dynamics}. We further consider the migration of two small thin films on the substrate curve, which is given by $y = 3 \sin \frac{2\pi}{15}x$. As shown in Fig.~\ref{fig:evolution_convexconcave}, in both cases we observe similar behaviours of the particles.

\begin{figure}[!htp]
  \centering
  \includegraphics[width=0.4\textwidth]{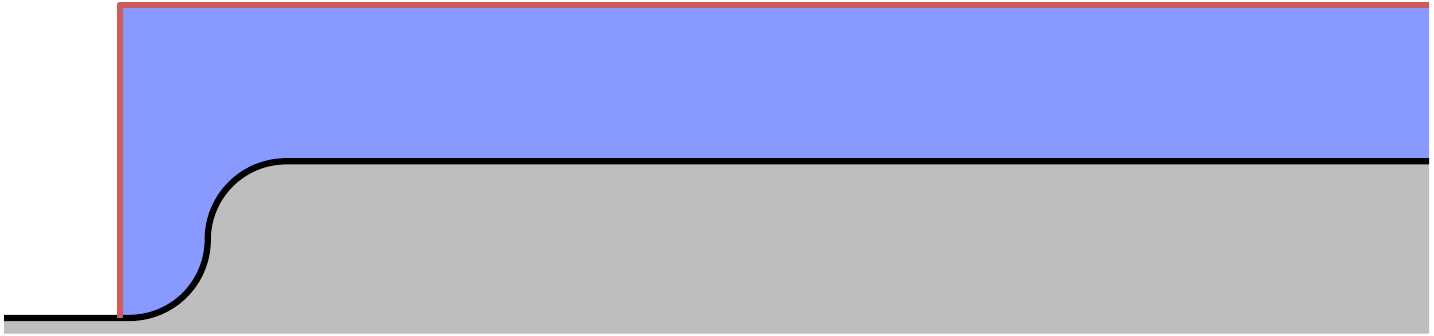}\hspace{0.9em}\includegraphics[width=0.4\textwidth]{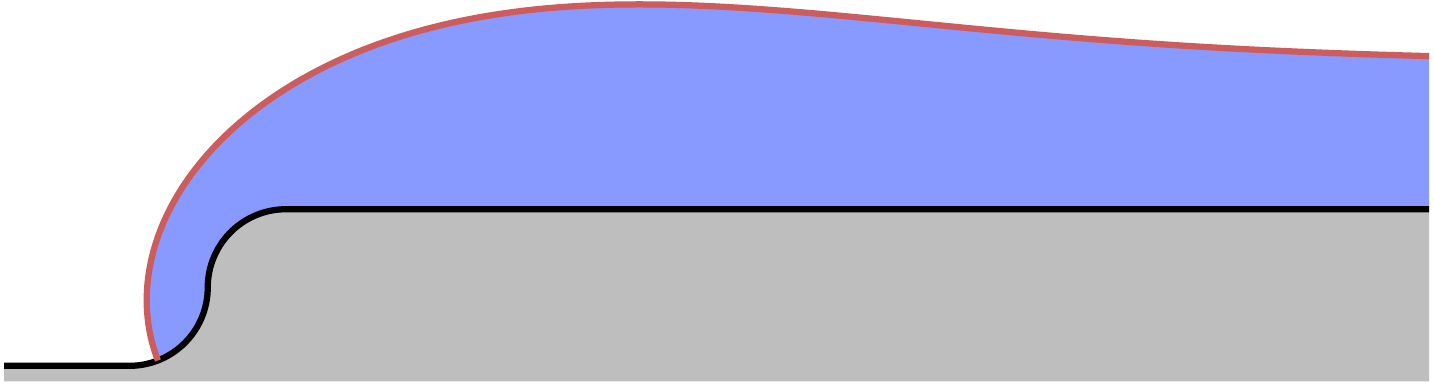}\\[0.8em]
  \includegraphics[width=0.4\textwidth]{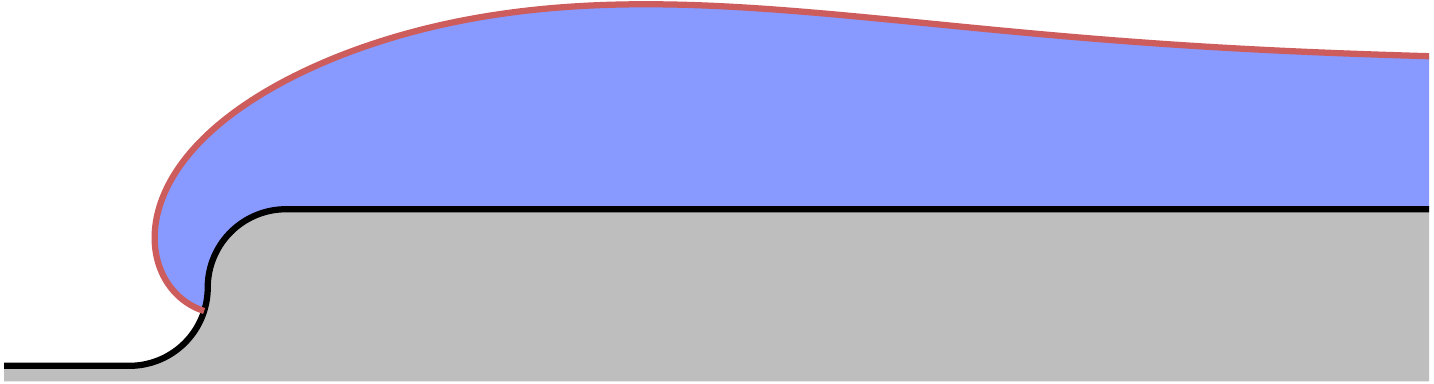}\hspace{0.9em}\includegraphics[width=0.4\textwidth]{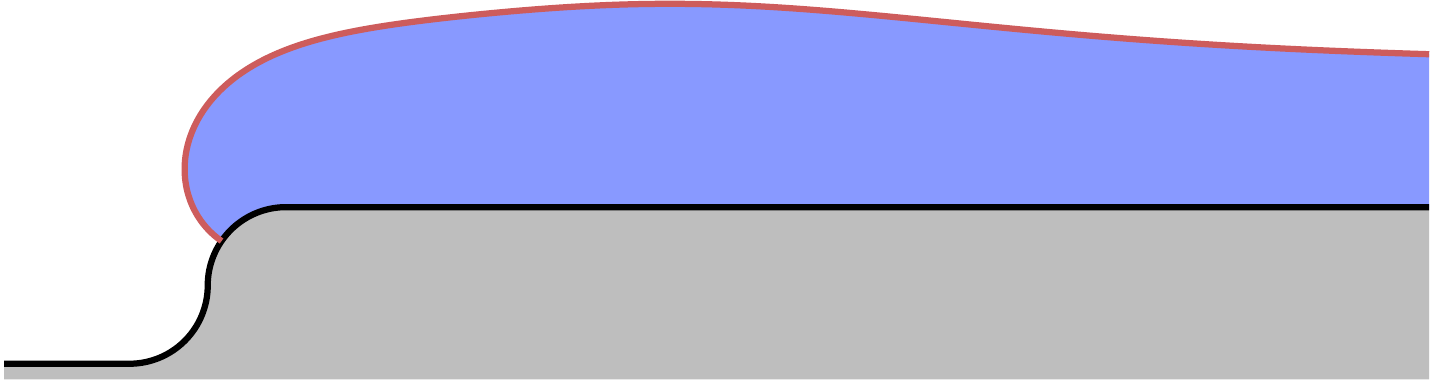}\\[0.8em]
  \includegraphics[width=0.4\textwidth]{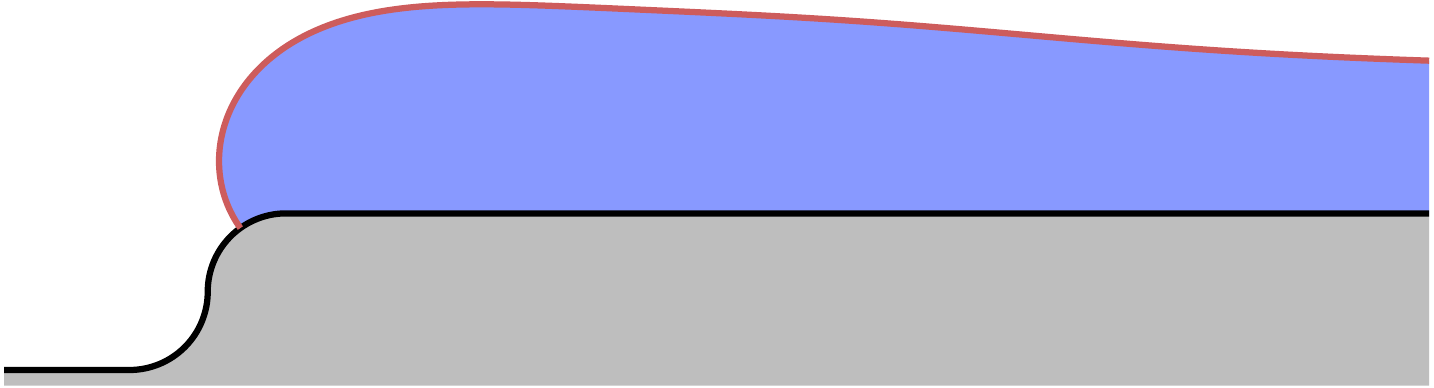}\hspace{0.9em}\includegraphics[width=0.4\textwidth]{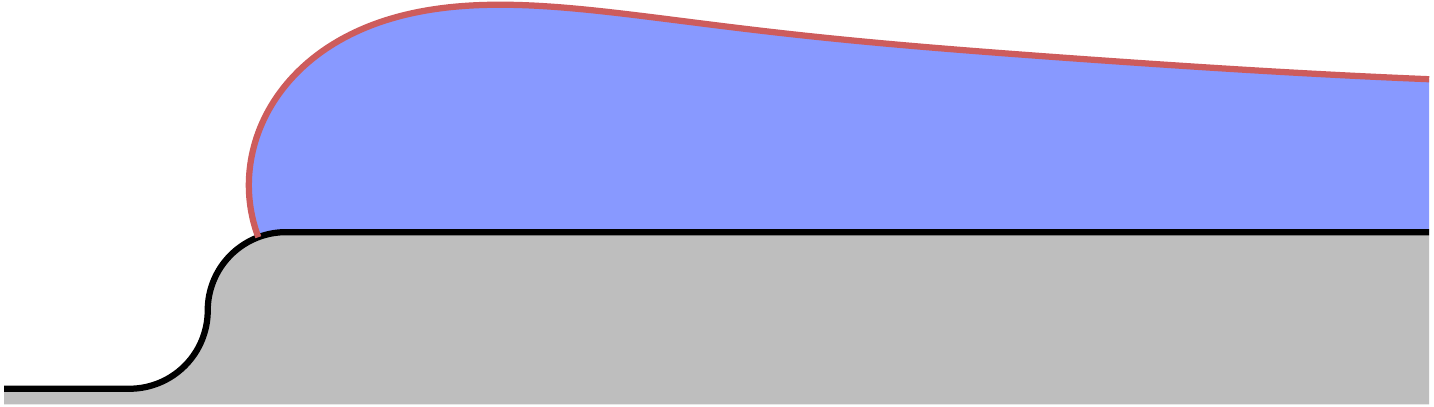}
  \caption{Snapshots in the edge migration of a semi-infinite step film across a corner at times $t = 0, 102.5, 105,  110,  122.5,  137.5 $ 
  (from left to right and then from top to bottom). }
  \label{fig:evolution_corner}
  \end{figure}

\begin{figure}[!htp]
\centering
\includegraphics[width=0.9\textwidth]{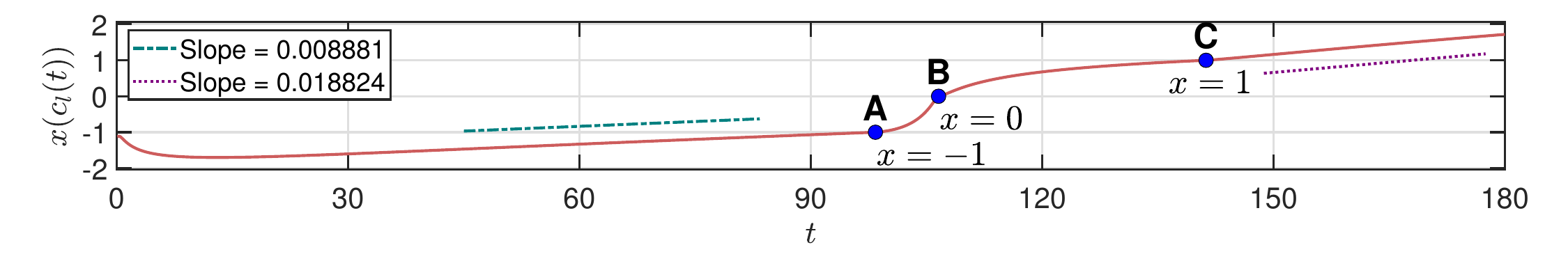}
\caption{The time history of the $x$-coordinate of the thin film edge, where points A, B, and C correspond to $x = -1, 0,$ and $1$. The green dash-dot and purple-dashed lines show the slope before A and the slope after C, respectively. }
\label{fig:evolution_corner_x}
\end{figure}

\vspace{0.4em}
\noindent
{\bf Example 5}: In our last example, we simulate the retraction of a semi-infinite step film across a corner. The initial step film is chosen as length $60$ and height $2$, positioned at a corner that has been smoothed by two circular arcs with radius $1$, as shown in  Fig. \ref{fig:evolution_corner}. Here we observe the edge retraction of the step film to approximate the corner and then gradually climbs to cross the corner. The evolution of x-coordinate of the thin film edge is also plotted in Fig.~\ref{fig:evolution_corner_x}, where we identify three special moments $A$, $B$, $C$. Here we observe an increase in the edge retraction rate after the thin film passed the corner, as indicated by the slopes in Fig. \ref{fig:evolution_corner_x}.  This phenomenon can be explained by the additional mass transport required to pass the corner, which is consistent to the results in \cite{Klinger12}.

\section{Conclusion}\label{sec:con}
We proposed a structure-preserving parametric finite element method for the sharp-interface model of solid-state dewetting with anisotropic surface energy on general curved substrates. The model is governed by anisotropic surface diffusion for an open curve with endpoints attached to a fixed smooth substrate curve. The introduced method is based on a novel weak formulation where the boundary conditions  can be naturally absorbed into the formulation or strongly enforced in the finite element spaces. The stability of the anisotropic surface energies is achieved using a time splitting of the Cahn-Hoffman vector with the help of a symmetrized surface energy matrix $\mathbf{Z}_k(\vec{n})$. We also adopted an arclength parameterization for the substrate curve, which helps to preserve the substrate energy stability on discrete level. The introduced method leads to a nonlinear system, and it can be solved by a hybrid iterative algorithm efficiently.  Moreover, our introduced method allows an tangential velocity to improve the mesh quality of the polygonal curve, and no remeshing procedure is required in practice.

\section*{Acknowledgement}
This work is partially supported by the Ministry of Education of Singapore under its AcRF Tier 2 funding MOE-T2EP20122-0002 (A-8000962-00-00) (W. Bao and Y. Li), and by the National Natural Science 
Foundation of China No. 12401572 (Q. Zhao). 

\bibliographystyle{model1b-num-names}
\bibliography{thebib}
\end{document}